\documentclass[notitlepage,leqno,10pt]{article}

\usepackage{latexsym}
\usepackage{amsmath}
\usepackage{amsthm}
\usepackage{amsfonts}
\usepackage{amssymb}
\usepackage{graphicx}
\renewcommand{\a }{\alpha }

\newcommand{\ve}{\varepsilon}

\newcommand{\rh }{\rho }

\newcommand{\intbar}{\mathop{\int\makebox(-13.5,0){\rule[4pt]{.7em}{0.3pt}}%
\kern-6pt}\nolimits}

\newcommand{\be}{\begin{equation}}
\newcommand{\ee}{\end{equation}}
\newenvironment{pf}{\noindent{\sc Proof}.\enspace}{\rule{2mm}{2mm}\medskip}

\newcommand{\R}{\mathbb{R}}
\newcommand{\Rn}{\mathbb{R}^n}

\newcommand{\Rdue}{\mathbb{R}^2}

\newcommand{\Rp}{\mathbb{R}^p}

\newcommand{\N}{\mathbb{N}}
\newcommand{\Sp}{\mathbb{S}}

\DeclareMathOperator{\diam}{diam}

\DeclareMathOperator{\graph}{graph}
\DeclareMathOperator{\C}{C}
\DeclareMathOperator{\D}{D}
\DeclareMathOperator{\de}{d}

\DeclareMathOperator{\dist}{dist}

\DeclareMathOperator{\Leins}{\mathcal{L}^1}
\DeclareMathOperator{\Lzwei}{\mathcal{L}^2}

\DeclareMathOperator{\B}{\mathcal{B}_\varepsilon}

\DeclareMathOperator{\diver}{div}

\DeclareMathOperator{\spt}{spt}

\begin{document}

\author{{\sc Ernst Kuwert, Andrea Mondino, Johannes Schygulla}}

\date{}

\title{Existence of immersed spheres minimizing curvature 
functionals in compact  $3$-manifolds}

\newtheorem{lem}{Lemma}[section]
\newtheorem{pro}[lem]{Proposition}
\newtheorem{thm}[lem]{Theorem}
\newtheorem{rem}[lem]{Remark}
\newtheorem{cor}[lem]{Corollary}
\newtheorem{df}[lem]{Definition}
\newtheorem*{Theorem}{Theorem}
\newtheorem*{Lemma}{Lemma}
\newtheorem*{Proposition}{Proposition}
\newtheorem*{claim}{Claim}

\maketitle

%\footnotetext[1]{Universit\"at Freiburg, mathematisches ??????? E-mail address: ????}
%\footnotetext[2]{SISSA, ?????? 34014 Trieste, Italy, E-mail address: mondino@sissa.it}

\begin{abstract}
We study curvature functionals for immersed $2$-spheres 
in a compact, three-dimensional Riemannian manifold $M$. Under the assumption that 
the sectional curvature $K^M$ is strictly positive, we prove the 
existence of a smooth immersion $f:\Sp^2 \to M$ minimizing the $L^2$ integral of 
the second fundamental form. Assuming instead that $K^M \leq 2$ and that 
there is some point $\overline{x} \in M$ with scalar curvature $R^M(\overline{x}) > 6$, 
we obtain a smooth minimizer $f:\Sp^2 \to M$ for the functional 
$\int \frac{1}{4}|H|^2+1$, where $H$ is the mean curvature.
\end{abstract}
\bigskip\bigskip

\begin{center}
\bigskip

\centerline{\bf AMS subject classification: }
53C21, 53C42, 58E99, 35J60
\end{center}

\section{Introduction}\label{s:in}
Let $M$ be a three-dimensional, compact Riemannian manifold with metric $h$. 
For any immersed closed surface $f:\Sigma \hookrightarrow M$ with induced metric 
$g = f^\ast h$ and second fundamental form $A$, we consider the functional 
\begin{equation} \label{eq:WE}
E(f) = \frac{1}{2} \int_\Sigma |A|^2\,d\mu_g.
\end{equation}
We denote by $H$ the mean curvature vector and by 
$A^\circ$ the tracefree component of $A$. The extrinsic curvature is related to 
the intrinsic curvature, i.e. the sectional curvature $K_g$ of the induced metric 
and the sectional curvature $K^M_f$ of the tangent plane in $TM$, by the Gau{\ss} 
equation 
\begin{equation} \label{eq:G}
\frac{1}{4}|H|^2 - \frac{1}{2}|A^\circ|^2
= \frac{1}{2} \big(|H|^2 - |A|^2\big) = K_g - K^M_f.
\end{equation}
Integrating and using the Gau{\ss}-Bonnet theorem yields the well-known identities 
\begin{equation} \label{eq:EQF}
\frac{1}{4} \int_\Sigma |H|^2\,d\mu_g + \frac{1}{2} \int_\Sigma |A^\circ|^2\,d\mu_g
= \frac{1}{2} \int_\Sigma |A|^2\,d\mu_g 
= \frac{1}{2} \int_\Sigma |H|^2\,d\mu_g + \int_\Sigma K^M_f\,d\mu_g - 2\pi \chi(\Sigma),
\end{equation}
where $\chi(\Sigma)$  is the Euler characteristic. For $M = \R^3$ the 
functional $E$ reduces to the classical Willmore energy given by 
\begin{equation} \label{eq:W}
W(f) = \frac{1}{4} \int_\Sigma |H|^2\,d\mu_g,
\end{equation}
more precisely we have $E(f) = 2 W(f) - 2\pi \chi(\Sigma)$. In \cite{Will} 
Willmore proved the inequality $W(f) \geq 4\pi$ for all $f:\Sigma \to \R^3$, with equality only for the round spheres.\\
%and conjectured that the torus 
%given by the equation $(r-\sqrt{2})^2 + z^2 = 1$ where $r = \sqrt{x^2+y^2}$ 
%is a minimizer.\\
\\
In the present paper we study the problem of minimizing $E(f)$ in the class
of immersed spheres in the Riemannian manifold $M$. Any totally geodesic 
$f:\Sp^2 \to M$ is trivially a minimizer, but totally geodesic immersions do 
not always exist. For instance, there are no totally umbilic surfaces 
in the Berger spheres (except $\Sp ^3$), see \cite{ST}. For appropriate parameters, these
spheres have positive sectional curvature \cite{Dan}. We prove the 
following existence result.

\begin{thm}\label{thm:ExEK}
Let $M$ be a compact, $3$-dimensional Riemannian manifold. On the class 
$[\Sp^2,M]$ of smooth immersions $f:\Sp^2 \to M$, consider the functional 
$$
E: [\Sp^2,M] \to \R,\,E(f) = \frac{1}{2} \int_{\Sp^2} |A|^2\,d\mu_g.
$$
If $M$ has sectional curvature $K^M > 0$, then there exists 
a minimizer $f$ in $[\Sp^2,M]$ for $E$. 
\end{thm}

We remark that our proof actually needs only the two conditions that 
$\inf_{[\Sp^2,M]} E(f) < 4\pi$ and that the area is bounded along some 
minimizing sequence. We always have $\inf_{[\Sp^2,M]} E(f) \leq 4\pi$, 
since the energy goes to $4\pi$ for a sequence of distance 
spheres shrinking to a point. Moreover, the strict inequality 
is necessary to rule out such a minimizing sequence. For example, 
if $M$ has strictly negative sectional curvature then $E(f) > 4\pi$ for 
any sphere immersed into $M$ by equation \eqref{eq:EQF}, and 
the infimum is not attained. Of course, the boundedness of 
the area along the minimizing sequence is also necessary, 
at least if we want subconvergence of the surface measures.
The first condition will be settled using a local expansion around
a point with strictly positive scalar curvature. The strong 
curvature assumption $K^M > 0$ of Theorem \ref{thm:ExEK}
is used to obtain the upper area bound. Possibly, the situation 
when the area actually goes to infinity (in the case when $K^M$ 
is not strictly positive) can be studied using results of 
Hutchinson \cite{Hu1} on curvature varifolds, see also \cite{MonVar}.\\ 
\\
In asymptotically flat $3$-manifolds $M$, spheres which are critical points 
of related curvature functionals have been constructed recently by Mondino \cite{Mon1,Mon2}
and Lamm, Metzger \& Schulze \cite{LMS}, see also \cite{LM}. They obtain the 
solutions as perturbations of round spheres using implicit function type arguments.
In \cite{SiL} L. Simon proved the existence of an embedded torus in $\R^n$, which 
minimizes the classical Willmore functional. Our approach implements his 
fundamental theory in the case of spheres immersed into the Riemannian manifold 
$M$. Recently, an alternative approach to Simon's theorem was developed
by Rivi\`{e}re \cite{Riv}.\\
\\
We now briefly outline the contents of the paper. In Section $2$, we gain some 
global control in terms of area and diameter bounds. For the lower diameter 
bound we use the bound $\inf E < 4\pi$ mentioned above. Local area 
bounds are then obtained by adapting Simon's monotonicity formula \cite{SiL}.
In Section $3$ we prove Theorem \ref{thm:ExEK}. First we obtain a limiting measure as a 
candidate for the minimizer. Adapting the arguments of \cite{SiL} to the Riemannian situation, 
we establish $C^{1,\alpha} \cap W^{2,2}$ regularity away from a finite
set of bad points where the curvature significantly concentrates. If a closed surface in $\R^3$ has Willmore energy below $8\pi$, 
as is the case in \cite{SiL}, then the area ratio is bounded below two 
by the monotonicity formula. Unfortunately, this involves a global argument 
which does not generalize immediately to our situation in $M$. We rule out the formation of branch points using the global bound
$\inf E < 4\pi$. This step involves a degree argument for the Gau{\ss} 
map, which does not extend to higher codimension. Eventually, we exclude all 
bad points and finally prove smoothness. The fact that
the limiting measure comes from an immersed sphere is proved using 
a compactness result of Breuning \cite {Breu}.\\
\\ 
In the final section $4$ we discuss the following variant of 
Theorem \ref{thm:ExEK}.
 
\begin{thm}\label{thm:ExW1}
For a closed, three-dimensional Riemannian manifold $M$, consider 
on the class of immersions $f:\Sp^2 \to M$ the functional 
$$
W_1(f) = \int_{\Sp^2} \Big(\frac{1}{4}|H|^2 +1\Big)\,d\mu_g.
$$
If $M$ has sectional curvature $K^M \leq 2$ and moreover the scalar curvature $R^M(\overline{x}) > 6$
for some point $\overline{x} \in M$, then there exists a smooth minimizer $f$ in $[\Sp^2,M]$ 
for $W_1$.
\end{thm}

We remark that the curvature conditions in Theorem \ref{thm:ExW1} can be 
fulfilled, for instance they hold for a round sphere $\Sp^3(R)$ if 
$\frac{1}{\sqrt{2}} \leq R < 1$. One motivation to study the functional is that
if we transform the classical Willmore integral from $\R^3$ to $\Sp^3$ using
stereographic projection, then we obtain the functional $W_1(f)$. Moreover, 
minimal surfaces are obvious critical points of $W_1$. The existence of 
minimizers for the functional $\int |H|^2\,d\mu_g$, possibly with branch 
points, was proposed in \cite{SiProc}. In our theorem, the assumption 
$K^M \leq 2$ in Theorem \ref{thm:ExW1} is mainly used to rule out the
branch points.
\begin{center}

{\bf Acknowledgments}

\end{center}

\noindent 
The second author would like to thank his supervisor Prof.\,A.\,Malchiodi for proposing to study the 
Willmore functional in Riemannian manifolds, and for his constant support. All authors acknowledge the
support by the DFG Collaborative Research Center SFB/Transregio 71 and of M.U.R.S.T, within the project 
B-IDEAS ''Analysis and Beyond'', making possible our cooperation with mutual visits at 
Freiburg and Trieste. 

\section{Global bounds for the minimizing sequence}\label{Sec:AprioriEst}

Here we collect some basic information for minimizing sequences of the functional $E$: 
global and local upper area bounds and a lower diameter bound. 
The first observation, following directly from \eqref{eq:EQF}, is 

\begin{pro}\label{prop:AreaEstEa1}
Let $M$ be a compact Riemannian $3$-manifold with sectional curvature $K^M >0$. 
%\begin{equation}\label{eq:HpMPos}
%\exists \lambda \text{ such that } \bar{K}>\lambda^2>0.
%\end{equation}
Then, for any immersed, closed surface $f: \Sigma \hookrightarrow M$ the 
total area $\mu_g(\Sigma)$ is bounded by 
\begin{equation}\label{eq:areaEst}
\mu_g(\Sigma) \leq C\,\Big(E(f) + 2\pi \chi(\Sigma)\Big) \quad \mbox{ with }
C= \frac{1}{\min_M K^M} < \infty.
\end{equation}
\end{pro}

We next apply Simon's monotonicity formula in $\R^m$ to show a local, quadratic 
area bound. 

\begin{lem}\label{lem:LocAreaEst}
Let $f:\Sigma \hookrightarrow M$ be a closed immersed surface in a 
compact $3$-manifold, with 
$$
W(f) + \mu_g(\Sigma) \leq \Lambda \quad \mbox{ for some } \Lambda < \infty.
$$
Then for any $x \in M$, $\rho > 0$ we have an estimate
$$
\mu_g\big(\{p \in \Sigma: f(p) \in B_\rho(x) \}\big) \leq C \rho^2, 
\quad \mbox{ where } C = C(\Lambda,M).
$$
\end{lem}

\begin{pf}
By Nash's theorem, there is an isometric embedding $I:M \hookrightarrow \R^m$ for some $m \in \N$.
The second fundamental forms of $f$, $I \circ f$ and $I$ are related by the formula
$$
A^{I \circ f}(\cdot,\cdot) = 
DI|_f \circ  A^{f}(\cdot,\cdot) \oplus (A^{I} \circ f) (Df,Df).
$$
Taking the trace and squaring yields for an orthonormal basis $v_i = Df \cdot e_i$
$$
|H^{I \circ f}|^2 = 
|H^f|^2 + \Big|\sum_{i = 1}^2 A^I \circ f (v_i,v_i)\Big|^2 
\leq |H^f|^2 + 2 |A^I|^2 \circ f.
$$
Integrating we see that 
$W(I \circ f) \leq W(f) + C\,\mu_g(\Sigma)$ where $C = \frac{1}{2} \max |A^I|^2$.
Thus for any $x \in M$, we get from Simon's monotonicity formula, see $(1.3)$ in \cite{SiL}, 
$$
\mu_g\big(\{p \in \Sigma: f(p) \in B^M_\rho(x)\}\big) \leq 
\mu_g\big(\{p \in \Sigma: I(f(p)) \in B_\rho^{\R^m}\big(I(x)\big)\}\big) 
\leq C \rho^2,
$$
with constant $C$ depending on $W(f)$, $\mu_g(\Sigma)$ and on $\max |A^I|$. \end{pf}

%Using the chain rule we compute at $p \in \Sigma$
%\begin{eqnarray*}
%\partial_j (\iota \circ f)(p) & = & (D\iota)\big(f(p)\big) \partial_j f(p),\\
%\partial^2_{ij} (\iota \circ f)(p) & = &
%(D^2 \iota)\big(f(p)\big) \big(\partial_i f(p),\partial_j f(p)\big) + (D\iota)\big(f(p)\big) D_i (\partial_j f)(p).
%\end{eqnarray*}
%Here we use the Levi-Civita connection $D$ on $M$ to define the second derivative, in particular we have 
%$D\iota \cdot D_X Y = D_X (D\iota \cdot Y)^{{\rm tan}_M}$. Now in Riemann normal coordinates at $p$, 
%we compute
%$$
%D_i (\partial_j) f = A^{f}_{ij} \quad \mbox{ and } \quad
%\partial^2_{ij} (\iota \circ f)(p) = A^{\iota \circ f}_{ij}(p). 
%$$
%Combining we arrive at
%$$
%A^{\iota \circ f}_{ij}(p) = A^{\iota}\big(f(p)\big) \big(\partial_i f(p),\partial_j f(p)\big) + (D\iota)\big(f(p)\big) A^{f}_{ij},
%$$
%and by taking traces
%$$
%\vec{H}^{\iota \circ f} = {\rm tr\,} A^{\iota}|_{{\rm im }Df(p)} \oplus D\iota\big(f(p)\big) \vec{H}^f.
%$$ 
Next we state an asymptotic expansion for the energy $E$ on geodesic spheres around 
a point $x \in M$, which follows from the well-known expansion of the metric in 
exponential coordinates. Since $|A|^2 = |A^\circ|^2 + \frac{1}{2} |H|^2$, we may 
combine Proposition $3.1$ in \cite{Mon1} with Lemma $3.5$ and Proposition $3.8$ in 
\cite{Mon2} to get the result. Note that for $M = \R^3$ we always have 
$E(f) \geq 4\pi$, with equality only for round spheres, by \cite{Will}.

\begin{lem}\label{lem:EaSpr}
Let $M$ be a $3$-dimensional Riemannian manifold. Then for geodesic spheres
$S_\rho(x) = \{y \in M: \dist(y,x) = \rho\}$ around $x \in M$ we have the expansion
$$
E\big(S_\rho(x)\big) = 4\pi - \frac{2\pi}{3} R^M(x) \rho^2 + \mathcal{O}(\rho^3) 
\quad \mbox{ as } \rho \searrow 0.
$$
In particular, if the scalar curvature $R^M(\overline{x}) > 0$ for some $\overline{x} \in M$, then 
$\inf_{f \in [\Sp^2,M]} E(f) < 4\pi$.
\end{lem}

At several points in this paper we work in local normal coordinates. 
The following lemma collects the relevant inequalities between the
Riemannian and the coordinate quantities. 

\begin{lem} \label{lem:localcomparison}
Let $h_{1,2}$ be Riemannian metrics on a manifold $M$, with norms satisfying
$$
(1+\ve)^{-1} \|\cdot \|_1 \leq \|\,\cdot \,\|_2 \leq (1+\varepsilon) \|\,\cdot \,\|_1
\quad \mbox{ for some } \ve \in (0,1].
$$
For any smooth immersed surface $f:\Sigma \to M$, the following 
inequalities hold with universal $C < \infty$: 
\begin{itemize}
\item[{\rm (i)}] ${\rm dist}_1(x,y) \leq (1+\ve) {\rm dist}_2(x,y)$ for all $x,y \in M$;
\item[{\rm (ii)}] $B^{h_1}_\sigma(x) \subset B^{h_2}_\rho(x)$,\, whenever $(1+\ve) \sigma \leq \rho$;
\item[{\rm (iii)}] $\mu_{g_1} \leq (1+C \ve) \mu_{g_2}$, where $g_{1,2} = f^\ast (h_{1,2})$; 
\item[{\rm (iv)}] $\|A_1\|_{1}^2 \leq \big(1+ C(\ve + \delta)\big) \|A_2\|_{2}^2 + C \delta^{-1} \|\Gamma\|^2_{h_1} \circ f$ 
                  for any $\delta \in (0,1]$, where $\Gamma := D^{h_1}-D^{h_2}$ and $D^{h_i}$ is the covariant derivative with respect to the metric $h_i$.
\end{itemize}
\end{lem}

\begin{pf} The statements (i) and (ii) are obvious. To compare the 
Jacobians of $f$ with respect to $h_{1,2}$, we use $\|\,\cdot \,\|_{g_1} \leq (1+\varepsilon) \|\,\cdot \,\|_{g_2}$
and compute for $v,w \in T_p \Sigma$ with $g_2(v,w) = 0$ 
$$
\|v \wedge w\|_{g_1}^2 = \|v\|_{g_1}^2 \|w\|_{g_1}^2 - g_1(v,w)^2 
\leq (1+\ve)^4  \|v\|_{g_2}^2 \|w\|_{g_2}^2 
= (1+\ve)^4  \|v \wedge w\|_{g_2}^2. 
$$
This proves the inequality (iii). Next we compare the norms for a bilinear map 
$B:T_p\Sigma \times T_p \Sigma \to T_{f(p)} M$. Choose a basis 
$v_\alpha$ of $T_p\Sigma$ such that $g_1(v_\alpha,v_\beta) = \delta_{\alpha \beta}$
and $g_2(v_\alpha,v_\beta) = \lambda_\alpha \delta_{\alpha \beta}$. 
Then
$$
\lambda_\alpha = \|v_\alpha\|_{g_2} \leq (1+ \ve) \|v_\alpha\|_{g_1} = 1 + \ve,
$$
and putting $w_\alpha = v_\alpha/\lambda_\alpha$ we obtain
$$
\|B\|_1^2 = \sum_{\alpha,\beta = 1}^2 
\lambda_\alpha^2 \lambda_\beta^2 \|B(w_\alpha,w_\beta)\|_{h_1}^2
\leq (1+C\ve) \sum_{\alpha,\beta = 1}^2 \|B(w_\alpha,w_\beta)\|_{h_2}^2
= (1+C\ve)\|B\|_2^2.
$$
Now denote by $P^\perp_{1,2}:T_{f(p)}M \to (T_p f)^{\perp_{h_{1,2}}}$ the orthogonal 
projections onto the normal spaces  with respect to $h_{1,2}$. Then we have for 
any $\delta > 0$ the estimate
\begin{eqnarray*}
\big\|A_1\big\|_1^2 & = & \big\|P_1^\perp D^{h_1}(Df)\big\|_1^2\\
& \leq & \big\|P_2^\perp D^{h_1}(Df)\big\|_1^2\\
& \leq & \big\|P_2^\perp \big(D^{h_2}(Df) + \Gamma\circ f (Df,Df)\big)\big\|_1^2\\
& \leq & (1+\delta) \big\|P_2^\perp D^{h_2}(Df)\big\|_1^2 + C \delta^{-1} \|\Gamma\|_{h_1}^2 \circ f\\
& \leq & (1+\delta)(1+C\varepsilon) \|A_2\|_2^2 + C \delta^{-1} \|\Gamma\|_{h_1}^2 \circ f.
\end{eqnarray*}
This proves the inequality (iv).
\end{pf}
                
The lower diameter bound follows by combining Proposition \ref{prop:AreaEstEa1},
Lemma \ref{lem:EaSpr} and the following fact. 
    
\begin{pro}\label{prop:LBdiamEa}
Let $M$ be a compact Riemannian $3$-manifold. Assume there is a 
minimizing sequence $f_k \in [\Sp^2,M]$ for $E(f)$ with $\diam f_k(\Sp^2) \to 0$
and $\mu_{g_k}(\Sigma) \leq C$. Then 
$$
\inf_{f \in [\Sp^2,M]} E(f) = 4\pi.
$$
\end{pro}

\begin{pf} After passing to a subsequence, we may assume that the
$f_k(\Sp^2)$ converge to a point $x_0 \in M$. For given $\varepsilon \in (0,1]$
we choose $\rho > 0$, such that in Riemann normal coordinates $x \in B_\rho(0) \subset \R^3$
$$
\frac{1}{1+\varepsilon} |\,\cdot \,| \leq \|\,\cdot \,\|_h \leq (1+\varepsilon) |\,\cdot \,| 
\quad \mbox{ and } \quad |\Gamma_{ij}^k(x)| \leq \varepsilon. 
$$
We have $f_k(\Sp^2) \subset B_\rho(x_0)$ for large $k$. Denoting by $A^e,g^e_k$ the 
quantities with respect to the coordinate metric, we get from Willmore's inequality
and Lemma \ref{lem:localcomparison}
$$
4\pi \leq \frac{1}{2} \int_\Sigma |A^e_{f_k}|_e^2\,d\mu_{g^e_k}
\leq (1+C\varepsilon) (1+\delta) \frac{1}{2} \int_\Sigma |A_{f_k}|^2\,d\mu_{g_k}
+ C(\delta) \varepsilon^2\, \mu_{g_k}(\Sigma).
$$
Since $\mu_{g_k}(\Sigma) \leq C$ by assumption, we may let first $k \to \infty$, 
then $\varepsilon \searrow 0$ and finally $\delta \searrow 0$ to obtain 
$$
\liminf_{k \to \infty} E(f_k) \geq 4\pi.
$$
As the upper bound follows from Lemma \ref{lem:EaSpr}, the lemma is proved. \end{pf}

\begin{lem} \label{lem:Link} Let $f:\Sigma \to M$ be a closed immersed 
surface in a compact $3$-manifold, and put $\Sigma_\rho(x_0) = f^{-1}(B_\rho(x_0))$ 
for $x_0 \in M$ and $\rho > 0$. There exist constants $\rho_0 > 0$ and $C < \infty$ 
depending only on $M$, such that for $x_0 \in f(\Sigma)$ we have 
$$
\frac{\mu_g(\Sigma_\sigma(x_0))}{\sigma^2} \leq 
C\,\Big(\frac{\mu_g(\Sigma_\rho(x_0))}{\rho^2} + \int_{\Sigma_\rho(x_0)} |H|^2\,d\mu_g\Big)
\quad \mbox{ whenever } 0 < \sigma \leq \rho \leq \rho_0.
$$
\end{lem}

\begin{pf} Again, we use an isometric embedding $I:M \to \R^m$. For $x_0 \in M$ 
we put
$$
\Sigma^{\R^m}_\rho(x_0) = (I \circ f)^{-1}\big(B^{\R^m}_\rho(I(x_0)).
$$
Choosing $\rho_0 > 0$ appropriately, we have 
$I(B_\rho(x_0)) \subset \left(B^{\R^m}_\rho(I(x_0))\cap I(M)\right)\subset I(B_{2\rho}(x_0))$ and 
hence
$$
\Sigma_\rho(x_0) \subset \Sigma_\rho^{\R^m}(x_0) \subset \Sigma_{2\rho}(x_0).
$$
Now from \cite{SiL}, we obtain for $0 < \sigma \leq \rho/2 \leq \rho_0$ the estimate
\begin{eqnarray*}
\frac{\mu_g(\Sigma_\sigma(x_0))}{\sigma^2} & \leq & \frac{\mu_g(\Sigma^{\R^m}_\sigma(x_0))}{\sigma^2}\\
& \leq & C\Big(\frac{\mu_g(\Sigma^{\R^m}_{\rho/2}(x_0))}{\rho^2} + 
\int_{\Sigma^{\R^m}_{\rho/2}(x_0)} |H^{I \circ f}|^2\,d\mu_g \Big)\\
& \leq &  C\frac{\mu_g(\Sigma_{\rho}(x_0))}{\rho^2} 
+ C \int_{\Sigma_{\rho}(x_0)} |H^f|^2\,d\mu_g
+ C \max |A^I|^2 \mu_g(\Sigma_{\rho}(x_0))\\
& \leq & C(1 + \rho_0^2 \max |A^I|^2) \frac{\mu_g(\Sigma_{\rho}(x_0))}{\rho^2}
+  C \int_{\Sigma_{\rho}(x_0)} |H^f|^2\,d\mu_g.
\end{eqnarray*}
This settles the inequality, if $\rho \geq 2\sigma$. As the claim is trivial 
for $\rho \in [\sigma,2\sigma]$, the lemma is proved. 
\end{pf}

\begin{lem} \label{lem:normalcoordinates}
Let $M$ be a Riemannian $3$-manifold, and $f:\Sigma \hookrightarrow M$ a closed
immersed surface with 
$$
W(f) + \mu_g(\Sigma) \leq \Lambda \quad \mbox{ for some } \Lambda < \infty.
$$
For any $\eta > 0$ there exist $\rho_0 = \rho_0(M,\eta) > 0$ and $C = C(M,\Lambda) < \infty$,
such that for any $x_0 \in M$, $x \in B_{\rho_0}(x_0)$ and $0 < \rho \leq \rho_0$ 
the following inequalities hold, where $B^e,g^e,\ldots$ are defined with respect to 
normal coordinates centered at $x_0$:
\begin{equation}
\label{eq:inclusion}
B_\sigma(x) \subset B^{e}_\rho(x),\, B^{e}_\sigma(x) \subset B_\rho(x) \quad
\mbox{ if } (1+ \eta)\sigma \leq \rho;
\end{equation}
\begin{equation}
\label{eqmeasurecomparison}
\frac{1}{1+\eta} \mu_{g^e}(\Sigma_\rho(x)) \leq \mu_g(\Sigma_\rho(x)) \leq (1+\eta) \mu_{g^e}(\Sigma_\rho(x));
\end{equation}
\begin{equation}
\label{eq:curvaturecomparison}
\frac{1}{1+\eta} \int_{\Sigma_\rho(x)} |A_e|_e^2\,d\mu_{g^e} - C \rho^2 \leq
\int_{\Sigma_\rho(x)} |A|^2\,d\mu_g \leq (1+\eta) \int_{\Sigma_\rho(x)} |A_e|_e^2\,d\mu_{g^e} + C \rho^2.
\end{equation}
%\item[{\rm (iii)}]
%\item[{\rm (iv)}]
\end{lem}

\begin{pf} We can assume that the assumption of Lemma \ref{lem:localcomparison} is 
satisfied on $B_{2\rho_0}(x_0)$ with $\ve = C(M) \rho_0^2$. The first two statements 
follow directly from that lemma. For (\ref{eq:curvaturecomparison}) we choose 
$\delta = \rho_0^2$ in Lemma \ref{lem:localcomparison}. Using 
$\|\Gamma\|_e \leq C \rho_0$, the statement follows by combining with Lemma 
\ref{lem:LocAreaEst}. \end{pf}

\section{Proof of Theorem \ref{thm:ExEK}}
For proving existence of a minimizer for the functional $E: [\Sp^2,M] \to \R,\,E(f) = \frac{1}{2} \int_{\Sp^2} |A|^2\,d\mu_g$, we use the direct method in the calculus of variations. Let $f_k:\Sp^2\hookrightarrow M$ be a minimizing sequence of immersed closed surfaces for the functional $E$. Denote by $\mu_k$ the Radon measure on $M$ given by
\begin{equation}
\mu_k(E) = \mu_{g_k}\big(f_k^{-1}(E)\big) = \int_E \theta_{f_k}(y)\,d{\cal H}^2(y),
\end{equation}
where $\theta_{f_k}$ is the multiplicity and $g_k$ is the induced metric.\\
\\
By Proposition \ref{prop:AreaEstEa1} we can assume
\begin{equation}
\mu_k \to \mu \quad \text{weakly as Radon measures}.
\end{equation}
Using this convergence and the monotonicity formula Lemma \ref{lem:Link}, it follows as in \cite{SiL} that
\begin{equation}
\spt \mu_k \to \spt \mu \quad \text{in the Hausdorff distance sense}.
\end{equation}
This Hausdorff convergence, together with Lemma \ref{lem:EaSpr} and Proposition \ref{prop:LBdiamEa}, implies that
\begin{equation}\label{lowdiam}
\diam_h (\spt \mu)\geq \liminf_k (\diam_h \spt \mu_k) >0.
\end{equation}
When working in normal coordinates, we denote the Euclidean coordinate quantities with an index $''e''$, for example $\mu_k^e$, $H_k^e$, $A_k^e$, $\ldots$, while the Riemannian quantities won't have any index.\\
\\
In order to prove regularity, we would like to apply Simon's Graphical Decomposition Lemma proved in \cite{SiL}. The most important assumption is that the $L^2$-norm of the second fundamental form is locally small, which we will need simultaneously for infinitely many $k$. Therefore we define the so called bad points with respect to a given $\varepsilon>0$ in the following way: Define the Radon measures $\alpha_k$ on $M$ by
$$\alpha_k=\mu_k \llcorner|A_k|^2.$$
Since $\alpha_k(M)\leq C$, there exists a Radon measure $\alpha$ on $M$ such that (after passing to a subsequence) $\alpha_k \to \alpha$ weakly as Radon measures. It follows that $\spt\alpha\subset \spt \mu$ and $\alpha(M)\leq C$. 

\begin{df}\label{df:BadPoints}
We define the bad points with respect to $\varepsilon>0$ by
\begin{equation*}
\B=\left\{\xi\in \spt \mu \,\big|\,\alpha(\{\xi\})>\varepsilon^2\right\}.
\end{equation*}
\end{df}

\begin{rem}\label{E7}
Since $\alpha(M)\le C$, there exist only finitely many bad points. Moreover if $\xi_0\in\spt\mu\setminus\B$, there exists a $\rho_0=\rho_0(\xi_0,\varepsilon)>0$ such that $\alpha(B_{\rho_0}(\xi_0))<2\varepsilon^2$, and since $\alpha_k  \to \alpha$ weakly we get
\begin{equation}\label{9}
\int_{B_{\rho_0}(\xi_0)} |A_k|^2\, d\mu_k \le 2\varepsilon^2\quad\text{for }k\text{ sufficiently large}.
\end{equation}
\end{rem}

From now on fix a point $\xi_0\in\spt\mu\setminus\B$ and choose normal coordinates around that point. In the following we will work in these fixed coordinates. Using the estimates in normal coordinates in Lemma \ref{lem:normalcoordinates} as well as Lemma \ref{lem:LocAreaEst}, the next Lemma is easily derived.

\begin{lem}\label{lem:EstGraDec}
For $\varepsilon\le\varepsilon_0$ there exists a $\rh_0=\rh_0(\xi_0, \varepsilon)>0$ and a $\beta=\beta(M)>0$, such that for all $\xi\in \spt \mu \cap B^e_{\frac{\rh_0}{2}}(\xi_0)$, for all $\rh \leq\frac{\rh_0}{4}$ and infinitely many $k$ 
\begin{eqnarray}
&i)& \mu^e_k(\overline{B^e_\rh(\xi)})\leq \beta \rh^2, \nonumber \\
&ii)& \int_{B^e_{\rh}(\xi)} |A^e_k|^2 d\mu^e_k \leq 3 \varepsilon^2. \nonumber
\end{eqnarray}   
\end{lem}

Thanks to Lemma \ref{lem:EstGraDec} we are in position to apply the Graphical Decomposition Lemma of Leon Simon (Lemma 2.1 in \cite{SiL}).

\begin{lem}\label{final}
For $\varepsilon\le\varepsilon_0$ there exists a $\rh_0=\rh_0(\xi_0,\varepsilon)>0$ such that for all $\xi\in \spt \mu \cap B_\frac{\rho_0}{2}^e(\xi_0)$, all $\rho\leq\frac{\rho_0}{4}$ and for infinitely many $k$ the following holds:
There exist 2-dimensional planes $L_l$, $1\le l\le M_k$, containing $\xi$; functions $u_k^l\in C^\infty(\overline{\Omega_k^l},L_l^\perp)$, where $\Omega_k^l=(B^e_\lambda(\xi)\cap L^l)\setminus\bigcup_m d_{k,m}^l$ with $\lambda>\frac{\rho}{4}$ and where the sets $d_{k,m}^l\subset L^l$ are pairwise disjoint closed discs disjoint from $\partial B^e_\lambda(\xi)$; sets $P_j^{k,l}\subset M$, $1\le j\le N_{k,l}$, which are diffeomorphic to discs and disjoint from $\graph u_k^l$; and open, connected sets $U_k^l\subset f_k^{-1}(B_{\frac{\rho}{4}}^e(\xi))$, such that
\begin{eqnarray*}
& (i) & D_k^l:=\graph u_k^l\cup\bigcup_{j=1}^{N_{k,l}}P_j^{k,l}\quad\text{is a topological disc}, \\
& (ii) & f_k(U_k^l)=D_k^l\cap B_{\frac{\rho}{4}}^e(\xi)=\Big(\graph u_k^l\cup\bigcup_{j=1}^{N_{k,l}}P_j^{k,l}\Big)\cap B_{\frac{\rho}{4}}^e(\xi), \\
& (iii) & f_k^{-1}(B_{\frac{\rho}{4}}^e(\xi))\text{ is the disjoint union of the sets }U_k^l.\phantom{\bigcup_{j=1}^{N_{k,l}}}
\end{eqnarray*}
%$$\mu_k^e\llcorner\overline{B_{\frac{\rho}{4}}^e(\xi)}=\sum_{l=1}^{M_k}\mathcal{H}^2_e\llcorner \Big(D_k^l\cap\overline{B_{\frac{\rho}{4}}^e(\xi)}\Big)=\sum_{l=1}^{M_k}\mathcal{H}^2_e\llcorner\Big(\Big(\graph u_k^l\cup\bigcup_{j=1}^{N_{k,l}}P_j^{k,l}\Big)\cap\overline{B_{\frac{\rho}{4}}^e(\xi)}\Big),$$
%where $\Omega_k^l=(B^e_\lambda(\xi)\cap L^l)\setminus\bigcup_m d_{k,m}^l$ with $\lambda>\frac{\rho}{4}$, and where the sets $d_{k,m}^l\subset L^l$ are pairwise disjoint closed discs disjoint from $\partial B^e_\lambda(\xi)$. Each set $D_k^l$ is a topological disc with $\graph u_k^l\cap \overline{B_{\frac{\rho}{4}}^e(\xi)} \subset D_k^l$ and $D_k^l\setminus \graph u_k^l$ is the union of the sets $P_j^{k,l}\subset f_k(\Sp^2)$, which are diffeomorphic to a closed disc. 
Moreover the following estimates hold:
\begin{equation}
M_k\le c=c(M),
\end{equation}
\begin{equation}
\sum_{l,m}\diam d_{k,m}^l+\sum_{l,j}\diam P_j^{k,l} \le c\left(\int_{B_{\rho}^e(\xi)}|A^e_k|^2\,d\mu^e_k\right)^\frac{1}{4}\rho\le c\varepsilon^\frac{1}{2}\rho,
\end{equation}
\begin{equation}
||u_k^l||_{L^\infty(\Omega_k^l)} \le c \varepsilon^{\frac{1}{6}}\rho+\delta_k,\quad||Du_k^l||_{L^\infty(\Omega_k^l)} \le c\varepsilon^{\frac{1}{6}}+\delta_k,\quad\text{where }\delta_k\to0.
\end{equation}
\end{lem}
Next we prove a lower $2$-density bound for the minimizing sequence $f_k$ away from the bad points, which we will need later. 
\begin{pro}\label{pro:LowerDensity}
For $\varepsilon\le\varepsilon_0$ there exists a $\rh_0=\rh_0(\xi_0, \varepsilon)>0$ and a constant $C=C(M)>0$, such that for all $\xi\in \spt \mu \cap B_{\rh_0}(\xi_0)$ and all $\rh \leq\rh_0$
$$\frac{\mu(B_\rh(\xi))}{\rh^2}\geq C .$$
\end{pro} 

\begin{pf}
Let $\rho_0=\rho_0(\xi_0,\varepsilon)>0$ as in Remark \ref{E7} and $\xi\in B_\frac{{\rho_0}}{2}(\xi_0)$. It follows that $B_{\frac{{\rho_0}}{2}}(\xi)\subset B_{\rho_0}(\xi_0)$. Choose according to the Hausdorff distance sense convergence a sequence $\xi_k\in \spt\mu_k$ such that $\xi_k\to \xi$. Therefore for given $\rho\le\rho_0$ and $k$ sufficiently large it follows that $B_\frac{\rho}{4}(\xi_k)\subset  B_\frac{\rho}{2}(\xi)\subset B_{\rh_0}(\xi_0)$. Since the norm of the mean curvature can be estimated by the norm of the second fundamental form, we get from (\ref{9}) for $k$ sufficiently large
\begin{equation*}
\int_{B_\frac{\rho}{4}(\xi_k)}|H_k|^2\,d\mu_k\leq c\int_{B_{\rh_0}(\xi_0)}|A_k|^2\,d\mu_k \leq c\varepsilon^2.
\end{equation*}
By letting $\sigma\to0$ in Lemma \ref{lem:Link}, it follows that
$$ 1 \leq  C\left( \frac{\mu_k(B_\frac{\rho}{4}(\xi_k))}{\rh^2} + \varepsilon^2 \right).$$
Choosing $\varepsilon_0^2\leq \frac{1}{2C}$ we get for $k$ sufficiently large $\frac{\mu_k(B_\frac{\rh}{2}(\xi))}{\rh^2}\geq C>0$, and the rest of the Proposition follows from the weak convergence $\mu_k\to\mu$.
%Now since $\mu$ is a finite Radon measure, for almost every $0<\rh\leq \rh_0$ we have $\mu(\partial B^g_{2\rh}(\xi))=0$ then  the weak convergence of measures implies $$\mu(B_{2\rh}^g(\xi))=\lim_k[\mu_k(B_{2\rh}^g(\xi))]=\lim_{k'}[\mu_{k'}(B_{2\rh}^g(\xi))].$$
%Since $\xi_{k'}\to \xi$, for $k'$ large enough $\xi_{k'}\in B_{\rh}^g(\xi)$ and $B_{\rh}^g(\xi_{k'})\subset B_{2 \rh}^g(\xi)$; it follows that 
%$$\lim_{k'}[\mu_{k'}(B_{2\rh}^g(\xi))]\geq \limsup_{k'}[\mu_{k'}(B_{\rh}^g(\xi_{k'}))]=\limsup_{k'}|\Sig_{k'}\cap B_{\rh}^g(\xi_{k'}) |_g\geq C\rh^2>0$$
%where in the last passage we used inequality \eqref{eq:LowDenSk}. Collecting the last two chains of inequalities we get
%$$\frac{\mu(B_{2\rh}^g(\xi))}{4\rh^2}\geq C>0$$
%for almost every $0<\rh \leq \rh_0$. Now fix an arbitrary $\rh\in (0,2 \rh_0)$, then there exists a sequence $\rh_n\uparrow \rh$  such that the last inequality is satisfied:  $\mu(B_{\rh_n}^g(\xi))\geq C \rh_n^2$. Passing to the limit in $n$ we get 
%$$\frac{\mu(B_{\rh}^g(\xi))}{\rh^2}\geq C>0 \quad \forall \; \rh\in (0,2 \rh_0). $$
%We can conclude using statement $i)$ of Proposition \ref{pro:Estge} and the smallness of $\rh_0$; indeed, for small $\rh$ we have $B_{\rh}^e(\xi)\supset B^g_{\rh+O(\rh^3)}(\xi)$ then 
%$$\mu(B^e_{\rh}(\xi)) \geq \mu(B^g_{\rh+O(\rh^3)}(\xi)) \geq C [\rh^2+O(\rh^6)] \geq C \rh^2 $$
% for all  $\rh\in (0,\rh_0)$, $\rh_0$ small enough. 
\end{pf}

%\begin{rem}\label{rem:muNot0}
%Proposition \ref{pro:LowerDensity} is crucial for our minimization problems since it avoids the trivial case when the candidate minimizer limit measure $\mu$ is null. Indeed in the case $\{f_k\}_{k\in \N}$ is a minimizing sequence for either $E$ or $W_1$, we know from Theorem \ref {thm:CompLscEa} ( respectively Theorem \ref{thm:CompLscW1}) that the support $\spt \mu$ of the limit measure $\mu$ is compact, connected and with positive diameter; hence it contains infinitely many points. Since for both the functionals the $L^2$ norms of the second fundamental forms of $f_k$ are uniformly bounded (for $E$ it is trivial, see Lemma \ref{lem:AboundedWe} for $W_1$), by Remark \ref{E7}, for every $\varepsilon >0$ there are infinitely many $\varepsilon$-good points. Thus, applying Proposition \ref{pro:LowerDensity}, we have that there exists a small $\rh_0>0$ such that $\mu(B^e_{\rh_0}(\xi_0))>C \rh_0^2>0$.   
%\end{rem}

In the next step we estimate the $L^2$-norm of the second fundamental form on small balls around the "good points". This estimate will help us to show that the candidate minimizer $\mu$ is actually the measure associated to $C^{1,\alpha}\cap W^{2,2}$-graphs in a neighborhood around the good points.
\begin{lem}\label{2ff-absch}
For $\varepsilon\leq \varepsilon_0$ there exists a $\rh_0=\rh_0(\xi_0,\varepsilon)>0$ such that for all $\xi\in\spt \mu \cap B_{\frac{\rho_0}{2}}^e(\xi_0)$ and all $\rho\le\frac{\rho_0}{4}$
$$\liminf_{k\to\infty}\int_{B_{\frac{\rho}{8}}^e(\xi)}|A^e_k|^2\,d\mu^e_k\le c\rho^\alpha,$$
where $c<\infty$ and $\alpha\in(0,1)$ only depend on the manifold $M$.
\end{lem}

\begin{pf} Let $\varepsilon\le\varepsilon_0$ such that Lemma \ref{lem:EstGraDec} and Lemma \ref{final} hold. Let $\rho_0=\rho_0(\xi_0,\varepsilon)>0$ as before  and apply the Graphical Decomposition Lemma for $\rho\le\frac{\rho_0}{4}$ given by Lemma \ref{final} to infinitely many $k$. For these $k$ (surface index),  $l\in\{1,\ldots,M_k\}$ (slice index) and $\gamma\in\left(\frac{\rho}{16},\frac{3\rho}{32}\right)$ define the set
$$C_{\gamma}^l(\xi)=\left\{x+y\,\big|\,x\in B^e_\gamma(\xi)\cap L_l, y\in L_l^\perp\right\}.$$
From the estimates for the diameters of the pimples $P_j^{k,l}$ and the $C^{1}$-estimates for the graph functions $u^l_k$, it follows that 
\begin{equation}\label{eq:DCB}
D_k^l\cap C_{\gamma}^l(\xi)=D_k^l\cap C_{\gamma}^l(\xi)\cap\overline{B_{\frac{\rho}{4}}^e(\xi)}\quad\text{for }\varepsilon\le\varepsilon_0\text{ and }\delta_k\text{ sufficiently small}.
\end{equation}
%To see this, let $z\in D_k^l\cap C_{\gamma}^l(\xi)$, then $z=x_1+y_1$ with $x_1\in B_\gamma(\xi)\cap L_l$, $y_1\in L_l^\perp$. Since $D^l_k$ is disjoint union of a graph and a pimple part, there are two possible cases:
%\begin{itemize}
%\item[1)] $z\in\graph u_k^l\cap C_{\gamma}^l(\xi)$: thus $|y_1|=|u_k^l(x_1)|\le c\varepsilon^\frac{1}{6}\rho+\delta_k$ and
%\begin{eqnarray}
%|z-\xi| & \le & |x_1-\xi|+|y_1| \phantom{\frac{3\rho}{4}} \le  \gamma+c\varepsilon^\frac{1}{6}\rho+\delta_k \phantom{\frac{3\rho}{4}} \le  \frac{3\rho}{32}+c\varepsilon^\frac{1}{6}\rho+\delta_k\phantom{\frac{3\rho}{4}} \nonumber \\
  %& \le & \frac{\rho}{8}+\delta_k\quad\text{for }\varepsilon\le\varepsilon_0, \; \varepsilon_0 \text{ maybe smaller } \nonumber \\
  %& \le & \frac{\rho}{4}\quad\text{for }\delta_k\le\frac{\rho}{8}. \label{eq:z-xi}
%\end{eqnarray}
%\item[2)] $z\in D_k^l\cap P_j^k\cap C_{\gamma}^l(\xi)$ for some $j\in\N$: Since $\diam P_j^k\le c\varepsilon^\frac{1}{2}\rho$ it follows that $|y_1|\le c\varepsilon^\frac{1}{6}\rho+\delta_k+\diam P_j^k\le c\varepsilon^\frac{1}{6}\rho+\delta_k$. Now the claim follows in the same way as above in 1).
%\end{itemize}  
Next define the set $A_k^l$ by
$$A_k^l(\xi)=\left\{\gamma\in\left(\frac{\rho}{16},\frac{3\rho}{32}\right)\,\Big|\,\partial C_{\gamma}^l(\xi)\cap\bigcup_j P_j^{k,l}=\emptyset\right\}.$$
For $\varepsilon\le\varepsilon_0$ it follows that
$$\Leins(A_k^l(\xi))\ge\frac{\rho}{32}-\sum_j\diam P_j^{k,l}\ge\frac{\rho}{32}-c\varepsilon^\frac{1}{2}\rho\ge\frac{\rho}{64}.$$
From Lemma \ref{selection} it follows that there exists a set $T_l\subset\left(\frac{\rho}{16},\frac{3\rho}{32}\right)$ with $\Leins(T_l)\ge\frac{\rho}{64}$, such that for all $\gamma\in T_l$ 
$$\partial C_{\gamma}^l(\xi)\cap\bigcup_j P_j^{k,l}=\emptyset\quad\text{for infinitely many }k.$$
Now let $\gamma\in T_l$ be arbitrary (it will be chosen later). We apply the Extension Lemma \ref{extension} given in the Appendix to get a function $w_k^l\in C^\infty\left(\overline{B^e_\gamma(\xi)}\cap L_l, L_l^\perp\right)$ such that
\begin{eqnarray*}
w_k^l  =  u_k^l & , & \frac{\partial w_k^l}{\partial\nu} = \frac{\partial u_k^l}{\partial\nu} \quad \text{on }\partial B^e_\gamma(\xi)\cap L_l,\phantom{\int_{d_{k,m}^\sim}}\\
||w_k^l||_{L^\infty(B^e_\gamma(\xi)\cap L_l)} & \le & c\varepsilon^\frac{1}{6}\gamma+\delta_k,\quad\text{where }\delta_k\to0,\phantom{\int_{d_{k,m}^\sim}}\\
||\D w_k^l||_{L^\infty(B^e_\gamma(\xi)\cap L_l)} & \le &  c\varepsilon^{\frac{1}{6}}+\delta_k,\quad\text{where }\delta_k\to0,\phantom{\int_{d_{k,m}^\sim}}\\
\int_{B^e_\gamma(\xi)\cap L_l}|\D^2 w_k^l|^2 & \le & c\gamma\int_{\graph {u_k^l}_{|_{\partial B^e_\gamma(\xi)\cap L_l}}}|A^e_k|^2\,d{\cal H}^1_e,
\end{eqnarray*}
where $d{\cal H}^1_e$ is the 1-dimensional Euclidean Hausdorff measure.\\
\\
Observe that, with an analogous argument as above using the estimates on $w^l_k$, we get
\begin{equation}\label{eq:wBe}
\graph w^l_k\subset\overline{B_{\frac{\rho}{4}}^e(\xi)}\quad\text{for }\varepsilon\le\varepsilon_0\text{ and }\delta_k\text{ sufficiently small}.
\end{equation}
By exchanging for each $l$ the disc $D_k^l\cap C_{\gamma}^l(\xi)$ with the disc $\graph w_k^l$, we get a new immersed surface $\tilde\Sigma_k$, which can be parametrized on $\Sp^2$ by a $C^{1,1}$-immersion $\tilde{f_k}:\Sp^2 \hookrightarrow M$. To simplify notation at this point and later in the paper we just write
\begin{equation}\label{def:tildeSigma}
\tilde{\Sigma}_k=\left(f_k(\Sp^2)\setminus\left(\bigcup_l D_k^l\cap C_{\gamma}^l(\xi)\right)\right)\cup\bigcup_l\graph w_k^l.
\end{equation}
%Now we consider the immersed surfaces 
%\begin{equation}\label{def:tildeSigma}
%\tilde{\Sigma}_k=\left(f_k(\Sp^2)\setminus\left(\bigcup_l D_k^l\cap C_{\gamma}^l(\xi)\right)\right)\cup\bigcup_l\graph w_k^l.
%\end{equation}
%Let us check that $\tilde{\Sigma}_k$ can be parametrized on $\Sp^2$ by a $C^{1,1}$-immersion $\tilde{f_k}:\Sp^2 \hookrightarrow M$: 
%Since the pimples $P_j^{k,l}$ are diffeomorphic to discs and since we have chosen a good radius $\gamma$ for the cylinder $C_\gamma ^l(\xi)$, 
%it is possible to show that $D_k^l\cap C_{\gamma}^l(\xi)$ is diffeomorphic to $\graph w_k^l$ for all $k,l$. Because of the boundary 
%properties of $w_k^l$, one can define a $C^{1,1}$-immersion $\tilde{f_k}:\Sp^2 \hookrightarrow M$ which parametrizes $\tilde{\Sigma}_k$. \\
More precisely we have to do the following: Choose a radius $\gamma'>\gamma$ such that the disc $D_k^l$
has a smooth graph representation by $u_k^l$ on the annulus $A:=B_{\gamma'}(\xi)\setminus\overline{B_\gamma(\xi)}\cap L_l$. 
Consider the disjoint union of the disc $B_{\gamma'}(\xi)\cap L_l$ and the topological disc 
$\Delta:=\Sp^2\setminus f_k^{-1}(D_k^l\cap\overline{C_{\gamma}^l(\xi)})$. Consider the diffeomorphism 
$\phi:A\to\Sp^2\setminus f_k^{-1}(D_k^l\cap C_{\gamma'}^l(\xi)\setminus\overline{C_{\gamma}^l(\xi)})$ 
given by $\phi(x)=f_k^{-1}(x+u_k^l(x))$. We define the smooth 2-manifold $\Sigma$ by identifying 
$x\in A$ and $p\in\Sp^2\setminus f_k^{-1}(D_k^l\cap C_{\gamma'}^l(\xi)\setminus\overline{C_{\gamma}^l(\xi)})$ 
if $\phi(x)=p$. We thus get a $C^{1,1}$-immersion $\tilde f_k:\Sigma \to M$ by putting
$$\tilde f_k= \begin{cases}
         \ f_k & \text{on } \Delta \\
         \ x+u_k^l(x) & \text{for } x\in A \\
         \ x+w_k^l(x) & \text{for } x\in\overline{B_\gamma(\xi)}\cap L_l.
     \end{cases}
$$     
It is easy to check that $\Sigma$ is orientable and has cohomology $H^1(\Sigma)=0$, and hence $\Sigma$ 
is diffeomorphic to $\Sp^2$. This constructs the desired $C^{1,1}$-immersion of $\Sp^2$.\\
\\
From the  definition of $\gamma$ we have that
\begin{equation*}
\int_{\graph w_k^l}|A_e|^2\,d{\cal H}_e^2  \le  c\int_{B^e_\gamma(\xi)\cap L_l}|\D^2 w_k^l|^2  \le c\gamma\int_{\graph {u_k^l}_{|_{\partial B^e_\gamma(\xi)\cap L_l}}}|A^e_k|^2 d{\cal H}_e^1 =  c\gamma\int_{\partial C_{\gamma}^l(\xi)\cap D_k^l}|A^e_k|^2d{\cal H}_e^1.
\end{equation*}
%Until now,  $\gamma\in T_l\subset\left(\frac{\rho}{16},\frac{3\rho}{32}\right)$ was arbitrary and $\Leins(T_l)\ge\frac{\rho}{64}$. Therefore, with a little Fubini-type argument,  we get  that the set 
%$$S_k^l=\left\{\gamma\in T_l\,\Big|\,\int_{\partial C_{\gamma}^l(\xi)\cap D_k^l}|\textrm{A}_e|^2 d\mu_e\le\frac{128}{\rho}\int_{\left(D_k^l\cap C_{k,\frac{3\rho}{32}}^l(\xi)\setminus C_{k,\frac{\rho}{16}}^l(\xi)\right)\setminus\bigcup_j P_j^k}|\textrm{A}_e|^2\,d\mu_e\right\}$$
%has measure $\Leins(S_k^l)\ge\frac{\rho}{128}$. Indeed otherwise we would have that
%\begin{eqnarray*}
%\int_{\left(D_k^l\cap C_{k,\frac{3\rho}{32}}^l(\xi)\setminus C_{k,\frac{\rho}{16}}^l(\xi)\right)\setminus\bigcup_j P_j^k}|\textrm{A}_e|^2\,d\mu_e & \ge & \int_{T_l\setminus S_k^l}\int_{\partial C_{\gamma}^l(\xi)\cap D_k^l}|\textrm{A}_e|^2 \\
%& \hspace{-5cm}\ge & \hspace{-2,5cm}\Leins\left(T_l\setminus S_k^l\right)\frac{128}{\rho}\int_{\left(D_k^l\cap C_{k,\frac{3\rho}{32}}^l(\xi)\setminus C_{k,\frac{\rho}{16}}^l(\xi)\right)\setminus\bigcup_j P_j^k}|\textrm{A}_e|^2\,d\mu_e \\
%& \hspace{-5cm}> & \hspace{-2,5cm}\left(\frac{\rho}{64}-\frac{\rho}{128}\right)\frac{128}{\rho}\int_{\left(D_k^l\cap C_{k,\frac{3\rho}{32}}^l(\xi)\setminus C_{k,\frac{\rho}{16}}^l(\xi)\right)\setminus\bigcup_j P_j^k}|\textrm{A}_e|^2\,d\mu_e \\
%& \hspace{-5cm}= & \hspace{-2,5cm}\int_{\left(D_k^l\cap C_{k,\frac{3\rho}{32}}^l(\xi)\setminus C_{k,\frac{\rho}{16}}^l(\xi)\right)\setminus\bigcup_j P_j^k}|\textrm{A}_e|^2\,d\mu_e,
%\end{eqnarray*}
%a contradiction. 
Until now $\gamma\in T_l\subset\left(\frac{\rho}{16},\frac{3\rho}{32}\right)$ was arbitrary. Since $\Leins(T_l)\ge\frac{\rho}{64}$, it easily follows from a simple Fubini-type argument as done in \cite{SiL} that we can choose $\gamma\in T_l$ such that for every $l,k$
$$\int_{\graph w_k^l}|A_e|^2\,d{\cal H}_e^2 \le c\int_{\left(D_k^l\cap C_{\frac{3\rho}{32}}^l(\xi)\setminus C_{\frac{\rho}{16}}^l(\xi)\right)\setminus\bigcup_j P_j^{k,l}}|A^e_k|^2\, d{\cal H}_e^2.$$
Now notice that for $\varepsilon\le\varepsilon_0$ (this follows from the estimates for $u_k^l$ and $\D u_k^l$)
\begin{eqnarray*}
B_{\frac{\rho}{16}}^e(\xi)\subset C_{\frac{\rho}{16}}^l(\xi)\quad\text{and}\quad\left(D_k^l\cap C_{\frac{3\rho}{32}}^l(\xi)\right)\setminus\bigcup_j P_j^{k,l}\subset\left(D_k^l\cap B_{\frac{\rho}{8}}^e(\xi)\right)\setminus\bigcup_j P_j^{k,l}.
\end{eqnarray*}
We get that
$$\int_{\graph w_k^l}|A_e|^2\,d{\cal H}_e^2\le c\int_{D_k^l\cap B_{\frac{\rho}{8}}^e(\xi)\setminus B_{\frac{\rho}{16}}^e(\xi)}|A^e_k|^2\,d{\cal H}_e^2.$$
By summing over $l$ and using the uniform bound on $M_k$ it follows that
\begin{equation}\label{eq:ProvEstAe}
\sum_{l=1}^{M_k}\int_{\graph w_k^l}|A_e|^2\,d{\cal H}_e^2 \le  c \sum_{l=1}^{M_k}\int_{D_k^l\cap B_{\frac{\rho}{8}}^e(\xi)\setminus B_{\frac{\rho}{16}}^e(\xi)}|A_k^e|^2\,d{\cal H}_e^2 = c\int_{B_{\frac{\rho}{8}}^e(\xi)\setminus B_{\frac{\rho}{16}}^e(\xi)} |A^e_k|^2\, d\mu^e_k. 
\end{equation}
Since $f_k$ is a minimizing sequence for the functional $E$ we get
$$E(\tilde {f_k})\ge E(f_k)-\varepsilon_k,\quad\text{where }\varepsilon_k\to0,$$
which implies
\begin{equation}\label{eq:S-Stilde}
\sum_{l=1}^{M_k} \int_{\graph w_k^l}|A|^2\,d{\cal H}^2\ge \int_{B^e_{\frac{\rh}{16}}(\xi)}|A_k|^2\,d\mu_k-\varepsilon_k.
\end{equation}
Using the estimates of Lemma \ref{lem:normalcoordinates} we finally get that
\begin{equation}\label{SS-Stilde}
\int_{B^e_{\frac{\rh}{16}}(\xi)} |A^e_k|^2\,d\mu^e_k \le c\int_{B_{\frac{\rho}{8}}^e(\xi)\setminus B_{\frac{\rho}{16}}^e(\xi)}|A^e_k|^2\,d\mu^e_k+\varepsilon_k+c\rho^2.
\end{equation}
By adding $c$ times the left hand side of this inequality to both sides ("hole-filling") we get that for all $\rho\le\frac{\rho_0}{4}$ and infinitely many $k$ 
$$\int_{B^e_{\frac{\rh}{16}}(\xi)} |A^e_k|^2\,d\mu^e_k \le \theta \int_{B_{\frac{\rho}{8}}^e(\xi)}|A^e_k|^2\,d\mu^e_k+\varepsilon_k+c\rho^2,$$
where $\theta=\frac{c}{c+1}\in(0,1)$ only depends on the manifold $M$. Now if we let $g(\rho)=\liminf_{k\to\infty}\int_{B_\rho^e(\xi)}|A^e_k|^2\,d\mu^e_k$, we get that
$$g(\rho)\le\theta g(2\rho)+c\rho^2\quad\text{for all }\rho\le\frac{\rho_0}{64}.$$
In view of Lemma \ref{decay} in the Appendix the present Lemma is proved.
\end{pf}

Now we are able to show that, in a neighborhood of the good points, the limit measure $\mu$ is the Radon measure associated to $C^{1,\alpha}\cap W^{2,2}$-graphs. First we recall the setting shortly: We had that $u_k^l:\Omega_k^l\to L_l^\perp$, where the set $\Omega_k^l$ was given by
$$\Omega_k^l=\left(B_\lambda^e(\xi)\cap L_l\right)\setminus\bigcup_m d_{k,m}^l,$$
where $\lambda>\frac{\rho}{4}$, and where the sets $d_{k,m}^l\subset L_l$ are pairwise disjoint closed discs which do not intersect $\partial B_\lambda^e(\xi)$.\\
\\
Define the quantity $\alpha_k(\rho)$ by
$$\alpha_k(\rho)=\int_{B^e_{4\rho}(\xi)}|A^e_k|^2\,d\mu^e_k,$$
and notice that by Lemma \ref{2ff-absch} and Lemma \ref{lem:EstGraDec} we have 
\begin{equation}\label{E29}
\liminf_{k\to\infty}\alpha_k(\rho)\le \min\left\{c\rho^\alpha,c\varepsilon^2\right\}\quad\text{for all }\rho\le\frac{\rho_0}{128}.
\end{equation}
Moreover it follows from Lemma \ref{final} that
\begin{equation}\label{E28}
\sum_m\diam d_{k,m}^l\le c\alpha_k(\rho)^\frac{1}{4}\rho.
\end{equation}
Therefore for $\varepsilon\le\varepsilon_0$ we may apply the generalized Poincar\'{e} inequality Lemma \ref{Poincare} to the functions $f_j^l=\D_j u_k^l$ and $\delta=c\alpha_k(\rho)^\frac{1}{4}\rho$ in order to get a constant vector $\eta_k^l$, with $|\eta_k^l|\le c\varepsilon^\frac{1}{6}+\delta_k\le c$ and $\delta_k\to0$, such that
$$\int_{\Omega_k^l}\left|\D u_k^l-\eta_k^l\right|^2\le c\rho^2\int_{\Omega_k^l}\left|\D^2 u_k^l\right|^2+c\alpha_k(\rho)^\frac{1}{4}\rho^2\sup_{\Omega_k^l}\left|\D u_k^l\right|^2.$$
Now we have
$$\int_{\Omega_k^l}\left|\D^2 u_k^l\right|^2\le c\int_{\graph u_k^l}|A^e_k|^2\,d{\cal H}_e^2 \le c\int_{B^e_{2\rho}(\xi)}|A^e_k|^2\,d\mu^e_k\le c\alpha_k(\rho).$$
Since $|\D u_k^l|\le c$ and $\alpha_k(\rho)\le1$ for $\varepsilon\le\varepsilon_0$, it follows that
\begin{equation}
\int_{\Omega_k^l}\left|\D u_k^l-\eta_k^l\right|^2\le c\alpha_k(\rho)^\frac{1}{4}\rho^2.
\end{equation}
Now let $\overline{u}_k^l\in C^{1,1}(B^e_\lambda(\xi)\cap L_l, L_l^\perp)$ be an extension of $u_k^l$ to all of $B^e_\lambda(\xi)\cap L_l$ as in Lemma  \ref{extension}, namely
\begin{eqnarray*}
\overline{u}_k^l&=&u_k^l\quad\text{in }B^e_\lambda(\xi)\cap L_l\setminus\bigcup_m d_{k,m}^l,\phantom{\frac{\partial\overline{u}_k^l}{\partial\nu}}\\
\overline{u}_k^l=u_k^l &\text{and}& \frac{\partial\overline{u}_k^l}{\partial\nu}=\frac{\partial u_k^l}{\partial\nu}\quad\text{on }\bigcup_m\partial d_{k,m}^l,\\
||\overline{u}_k^l||_{L^\infty(d_{k,m}^l)} & \le & c\varepsilon^\frac{1}{6}\rho+\delta_k,\quad\text{where }\delta_k\to0,\phantom{\frac{\partial\overline{u}_k^l}{\partial\nu}\bigcup_m d_{k,m}^l}\\   
||\D\overline{u}_k^l||_{L^\infty(d_{k,m}^l)} & \le &c\varepsilon^{\frac{1}{6}}+\delta_k,\quad\text{where }\delta_k\to0.\phantom{\frac{\partial\overline{u}_k^l}{\partial\nu}\bigcup_m d_{k,m}^l}
\end{eqnarray*}
It follows that $||\overline{u}_k^l||_{L^\infty(B_\lambda^e(\xi)\cap L_l)}+ ||\D\overline{u}_k^l||_{L^\infty(B_\lambda^e(\xi)\cap L_l)} \le c$, where $c$ is independent of $k$.
From the gradient estimates for the function $\overline{u}_k^l$, since $|\eta_k^l|\le c$ and because of (\ref{E28}) we get that
%\begin{eqnarray*}
%\int_{B_\lambda^e(\xi)\cap L_l}\left|\D\overline{u}_k^l-\eta_k^l\right|^2 & = & \int_{\Omega_k^l}\left|\D u_k^l-\eta_k^l\right|^2+\sum_m\int_{d_{k,m}^l}\left|\D\overline{u}_k^l-\eta_k^l\right|^2 \\
%& \le & c\alpha_k(\rho)^\frac{1}{4}\rho^2+c\sum_m\int_{d_{k,m}^l}\left|\D\overline{u}_k^l\right|^2+c\sum_m\int_{d_{k,m}^l}\left|\eta_k^l\right|^2\\
%& \le & c\alpha_k(\rho)^\frac{1}{4}\rho^2+c\sum_m\Lzwei(d_{k,m}^l)\le  c\alpha_k(\rho)^\frac{1}{4}\rho^2+c\left(\sum_m\diam d_{k,m}^l\right)^2 \\
%& \le & c\alpha_k(\rho)^\frac{1}{4}\rho^2+c\alpha_k(\rho)^\frac{1}{2}\rho^2 \le c\alpha_k(\rho)^\frac{1}{4}\rho^2.
%\end{eqnarray*}
\begin{equation}
\int_{B_\lambda^e(\xi)\cap L_l}\left|\D\overline{u}_k^l-\eta_k^l\right|^2\le c\alpha_k(\rho)^\frac{1}{4}\rho^2.
\end{equation}
Thus, in view of \eqref{E29}, we conclude that
\begin{equation}\label{E30}
\liminf_{k\to\infty}\int_{B_\lambda^e(\xi)\cap L_l}\left|\D\overline{u}_k^l-\eta_k^l\right|^2\le \min\left\{c\rho^{2+\alpha},c\varepsilon^\frac{1}{2}\rho^2\right\}\quad\text{for all }\rho\le\frac{\rho_0}{128}.
\end{equation}
Moreover, it trivially follows that $\left\|\overline{u}_k^l\right\|_{W^{1,2}(B_\lambda^e(\xi)\cap L_l)}\le c\rho^2\le c.$
Therefore it follows that the sequence $\overline{u}_k^l$ is equicontinuous and uniformly bounded in $C^1(B_\lambda^e(\xi)\cap L_l,L_l^\perp)$ and $W^{1,2}(B_\lambda^e(\xi)\cap L_l,L_l^\perp)$, and we get the existence of a function $u_\xi^l\in C^{0,1}(B_\lambda^e(\xi)\cap L_l,L_l^\perp)$ such that (after passing to a subsequence)
\begin{eqnarray*}
\overline{u}_k^l & \to & u_\xi^l\quad\text{in }C^0(B_\lambda^e(\xi)\cap L_l,L_l^\perp), \\
\overline{u}_k^l & \rightharpoonup & u_\xi^l\quad\text{weakly in }W^{1,2}(B_\lambda^e(\xi)\cap L_l,L_l^\perp),
\end{eqnarray*}
and such that the function $u_\xi^l$ satisfies the estimates
$$ \frac{1}{\rho} ||u_\xi^l||_{L^\infty(B_\lambda^e(\xi)\cap L_l)}+||\D u_\xi^l||_{L^\infty(B_\lambda^e(\xi)\cap L_l)}  \le  c\varepsilon^\frac{1}{6}.$$
Be aware that, a priori, the limit function might depend on the point $\xi$. Indeed, the sequence $u_k^l$ depends on $\xi$ since it comes from the Graphical Decomposition Lemma which is a local statement.\\
\\
Observe that, up to subsequences,  $\eta_k^l\to\eta^l$ with $|\eta^l|\le c\varepsilon^\frac{1}{6}$. Since $\D\overline{u}_k^l\rightharpoonup\D u_\xi^l$ weakly in $L^2(B_\lambda^e(\xi)\cap L_l)$ it follows that $\D\overline{u}_k^l-\eta_k^l\rightharpoonup\D u_\xi^l-\eta^l$ weakly in $L^2(B_\lambda^e(\xi)\cap L_l)$, and by lower-semicontinuity, estimate \eqref{E30} implies that
\begin{equation}\label{E31}
\int_{B_\lambda^e(\xi)\cap L_l}\left|\D u_\xi^l-\eta^l\right|^2\le \min\left\{c\rho^{2+\alpha},c\varepsilon^\frac{1}{2}\rho^2\right\}\quad\text{for all }\rho\le\frac{\rho_0}{128}.
\end{equation}
Now we can prove the graphical representation of the limit measure $\mu$.

\begin{lem}\label{mu=graph}
For $\varepsilon \leq \varepsilon_0$ there exists a $\rh_0=\rh_0(\xi_0,\varepsilon)$ such that for all $\xi\in\spt \mu\cap B_\frac{\rho_0}{2}^e(\xi_0)$ and all $\rho\le\rho_0$ we have
$$\mu\llcorner B_{\rho}^e(\xi)=\sum_{l=1}^M{\cal H}^2\llcorner\left(\graph u_\xi^l\cap B_{\rho}^e(\xi)\right),$$
where ${\cal H}^2$ denotes the $2$-dimensional Hausdorff measure of the Riemannian manifold $M$, and where each function $u_\xi^l\in C^{0,1}\left(B^e_{\rho}(\xi)\cap L_l,L_l^\perp\right)$ is as above, in particular
$$\frac{1}{\rho}||u_\xi^l||_{L^\infty(B_{\rho}^e(\xi)\cap L_l)}+\left\|\D u_\xi^l\right\|_{L^\infty(B^e_{\rho}(\xi)\cap L_l)}\le c \varepsilon^\frac{1}{6}.$$
\end{lem}

\begin{pf}
First we claim that for all $\rho\le\frac{\rho_0}{128}$ we have
\begin{equation}\label{E19}
\mu_k\llcorner B_\rho^e(\xi)=\sum_{l=1}^M{\cal H}^2\llcorner\left(\graph\overline{u}_k^l\cap B_\rho^e(\xi)\right)+\theta_k,
\end{equation}
where $\theta_k$ is a signed measure with $\liminf_{k\to\infty}$ of the total mass is smaller than $\min\left\{c\rho^{2+\alpha},c\varepsilon\rho^2\right\}$, i.e. $\theta_k=\theta_k^1-\theta_k^2$ with $\liminf_{k\to\infty}\left(\theta_k^1(M)+\theta_k^2(M)\right)\le \min\left\{c\rho^{2+\alpha},c\varepsilon\rho^2\right\}$.\\
\\
To prove the claim recall that the diameter estimates in Lemma \ref{final}, the quadratic area decay and the monotonicity formula Lemma \ref{lem:Link} yield 
$$\sum_{m,l}\Lzwei(d_{k,m}^l)+\sum_{j,l}{\cal H}^2(P_j^{k,l}) \le c\alpha_k(\rho)^\frac{1}{2}\rho^2.$$ Thus Lemma \ref{2ff-absch} yields for $\rho\le\frac{\rho_0}{128}$  $$\liminf_{k\to\infty}\sum_{m,l}\Lzwei(d_{k,m}^l)+\liminf_{k\to\infty}\sum_{j,l}{\cal H}^2(P_j^{k,l}) \le \min\left\{c\rho^{2+\alpha},c\varepsilon\rho^2\right\}.$$
The Graphical Decomposition Lemma \ref{final} yields $\mu_k\llcorner B_\rho^e(\xi)=\sum_{l=1}^M{\cal H}^2\llcorner\left(\graph\overline{u}_k^l\cap B_\rho^e(\xi)\right)+\theta_k$, where 
$$\theta_k=\sum_{l=1}^M{\cal H}^2\llcorner\left(\left(D_k^l\setminus\graph\overline{u}_k^l\right)\cap B_\rho^e(\xi)\right)-\sum_{l=1}^M{\cal H}^2\llcorner\left(\left(\graph\overline{u}_k^l\setminus D_k^l\right)\cap B_\rho^e(\xi)\right)=\theta_k^1-\theta_k^2.$$
We have that $\theta_k^1(M)\le\sum_{j,l}{\cal H}^2(P_j^{k,l})$ and $\theta_k^2(M)\le c\sum_{m,l}\Lzwei(d_{k,m}^l)$, and (\ref{E19}) follows.\\ 
\\
Now by taking limits in the measure theoretic sense we claim that
\begin{equation}\label{E20}
\mu\llcorner B_\rho^e(\xi)=\sum_{l=1}^M{\cal H}^2\llcorner\left(\graph u_\xi^l\cap B_\rho^e(\xi)\right)+\theta_\xi,
\end{equation}
where $\theta_\xi$ is a signed measure with total mass smaller than $\min\left\{c\rho^{2+\alpha},c\varepsilon^\frac{1}{4}\rho^2\right\}$. This equation holds for all $\rho\le\frac{\rho_0}{128}$ such that 
$$\mu\left(\partial B_\rho^e(\xi)\right)={\cal H}^2\llcorner\graph u_\xi^l\left(\partial B_\rho^e(\xi)\right)=0\quad\text{for all }l,$$
which holds for almost every $\rho$.\\
\\
To prove (\ref{E20}) let $U\subset M$ be an open subset.\\
\\
1.) Let $\rho\le\frac{\rho_0}{128}$ such that $\mu\left(\partial B^e_\rho(\xi)\right)=0$. Moreover assume that $\mu\llcorner B^e_\rho(\xi)\left(\partial U\right)=0$. It follows that $\mu\left(\partial\left(U\cap B^e_\rho(\xi)\right)\right)=0$ and therefore $\mu_k\left(U\cap B^e_\rho(\xi)\right)\to\mu\left(U\cap B^e_\rho(\xi)\right)$.\\
\\
2.) Let $\rho\le\frac{\rho_0}{128}$ such that ${\cal H}^2\llcorner\graph u_\xi^l\left(\partial B^e_\rho(\xi)\right)=0$. Assume that ${\cal H}^2\llcorner(\graph u_\xi^l\cap B^e_\rho(\xi))\left(\partial U\right)=0$. Now in general it follows for the $2$-dimensional Hausdorff measure of some $C^{0,1}$-graph $u$ that
$${\cal H}^2(\graph u)=\int\sqrt{\det g}=\int\sqrt{A(x,u(x))+B_i(x,u(x))\partial_i u(x)+C_{ij}(x,u(x))\partial_i u(x)\partial_j u(x)},$$
where the coefficients $A, B_i, C_{ij}$ just depend on the metric $h$ and are uniformly bounded in terms of the manifold $M$. Especially for the coefficient $A$ we have 
\begin{equation}
A(x)=h_{11}(x)h_{22}(x)-h_{12}(x)^2,
\end{equation}
where $h_{ij}$ are the coefficients of the metric $h$ of $M$. Therefore we get that the coefficient $A$ is bounded from below by a positive constant, namely there exists a constant $c_0>0$ such that
\begin{equation}\label{boundA}
\sup_{x\in M}A(x)\ge c_0>0.
\end{equation}
Using the $L^\infty$-bounds for the functions $\overline{u}_k^l$ and the coefficients $A, B_i, C_{ij}$, we get that 
\begin{eqnarray*}
& & \left|{\cal H}^2\llcorner\left(\graph\overline{u}_k^l\cap B^e_\rho(\xi)\right)(U)-{\cal H}^2\llcorner\left(\graph u_\xi^l\cap B^e_\rho(\xi)\right)(U)\right| \\
& & \hspace{0,5cm}\le c\int_{L_l}\left|\chi_{\phantom{}_{U\cap B^e_\rho(\xi)}}(x,\overline{u}_k^l(x))-\chi_{\phantom{}_{U\cap B^e_\rho(\xi)}}(x,u_\xi^l(x))\right|+\int_{L_l}\chi_{\phantom{}_{U\cap B^e_\rho(\xi)}}(x,u_\xi^l(x))\left|\sqrt{\det\overline{g}_k^l}-\sqrt{\det g^l}\right|
\end{eqnarray*}
Since $\overline{u}_k^l\to u_\xi^l$ uniformly and since ${\cal H}^2\llcorner\graph u_\xi^l\left(\partial B^e_\rho(\xi)\right)=0$, it follows that
$$\chi_{\phantom{}_{U\cap B^e_\rho(\xi)}}(x,\overline{u}_k^l(x))\to\chi_{\phantom{}_{U\cap B^e_\rho(\xi)}}(x,u_\xi^l(x))\quad\text{for a.e. }x\in L_l.$$
The Dominated Convergence Theorem yields 
$$\int_{L_l}\left|\chi_{\phantom{}_{U\cap B^e_\rho(\xi)}}(x,\overline{u}_k^l(x))-\chi_{\phantom{}_{U\cap B^e_\rho(\xi)}}(x,u_\xi^l(x))\right|\to0.$$
Now because of (\ref{boundA}) and the bounds for the functions $\overline{u}_k^l$ and $u_\xi^l$ it follows that for $\varepsilon\le\varepsilon_0$
\begin{eqnarray*}
\int_{L_l}\chi_{\phantom{}_{U\cap B^e_\rho(\xi)}}(x,u_\xi^l(x))\left|\sqrt{\det\overline{g}_k^l}-\sqrt{\det g^l}\right| & \le & \int_{L_l}\chi_{\phantom{}_{U\cap B^e_\rho(\xi)}}(x,u_\xi^l(x))\left|\det\overline{g}_k^l-\det g^l\right| \\
& \hspace{-11cm}\le & \hspace{-5,5cm}\int_{L_l}\chi_{\phantom{}_{U\cap B^e_\rho(\xi)}}(x,u_\xi^l(x))\left|A(x,\overline{u}_k^l(x))-A(x,u_\xi^l(x))\right| \\
& & \hspace{-5,5cm}+\int_{L_l}\chi_{\phantom{}_{U\cap B^e_\rho(\xi)}}(x,u_\xi^l(x))\left|B_i(x,\overline{u}_k^l(x))\partial_i\overline{u}_k^l(x)-B_i(x,u_\xi^l(x))\partial_i u_\xi^l(x)\right| \\
& & \hspace{-5,5cm}+\int_{L_l}\chi_{\phantom{}_{U\cap B^e_\rho(\xi)}}(x,u_\xi^l(x))\left|C_{ij}(x,\overline{u}_k^l(x))\partial_i\overline{u}_k^l(x)\partial_j\overline{u}_k^l(x)-C_{ij}(x,u_\xi^l(x))\partial_i u_\xi^l(x)\partial_j u_\xi^l(x)\right| \\
& \hspace{-11cm}=: & \hspace{-5,5cm}(1)+(2)+(3).\phantom{\int_{L_l}}
\end{eqnarray*}
Now $(1)\to0$ for $k\to\infty$ because of the uniform convergence $\overline{u}_k^l\to u_\xi^l$. The second term can be estimated in view of the boundedness of the coefficients $B_i$ and the functions $\overline{u}_k^l$ by
\begin{eqnarray*}
(2) & \le & c\int_{B^e_\rho(\xi)\cap L_l}\left|B_i(x,\overline{u}_k^l(x))-B_i(x,u_\xi^l(x))\right|+c\int_{L_l}\chi_{\phantom{}_{U\cap B^e_\rho(\xi)}}(x,u_\xi^l(x))\left|\D\overline{u}_k^l(x)-\D u_\xi^l(x)\right| \\
& \le & c\int_{B^e_\rho(\xi)\cap L_l}\left|B_i(x,\overline{u}_k^l(x))-B_i(x,u_\xi^l(x))\right|+c\int_{L_l}\chi_{\phantom{}_{U\cap B^e_\rho(\xi)}}(x,u_\xi^l(x))\left|\D\overline{u}_k^l(x)-\eta_k^l\right| \\
& & +c\int_{L_l}\chi_{\phantom{}_{U\cap B^e_\rho(\xi)}}(x,u_\xi^l(x))\left|\eta_k^l-\eta^l\right|+c\int_{L_l}\chi_{\phantom{}_{U\cap B^e_\rho(\xi)}}(x,u_\xi^l(x))\left|\eta^l-\D u_\xi^l(x)\right|
\end{eqnarray*}
The first term goes to 0, again by the uniform convergence $\overline{u}_k^l\to u_\xi^l$. For the second term we have that  $\chi_{\phantom{}_{U\cap B^e_\rho(\xi)}}(x,u_\xi^l(x))=0$ if $x\notin B^e_{\left(1-c\varepsilon^\frac{1}{6}\right)\rho}(\xi)\cap L_l$. This follows from the $L^\infty$-bound for the function $u_\xi^l$. Therefore we get that $$\left(\int_L\chi_{\phantom{}_{U\cap B^e_\rho(\xi)}}(x,u_\xi^l(x))\right)^\frac{1}{2}\le\Lzwei\left(B^e_{\left(1-c\varepsilon^\frac{1}{6}\right)\rho}(\xi)\cap L_l\right)^\frac{1}{2}\le c\rho.$$
In view of (\ref{E30}) we get $\liminf_{k\to\infty}\int_{L_l}\chi_{\phantom{}_{U\cap B^e_\rho(\xi)}}(x,u_\xi^l(x))|\D\overline{u}_k^l(x)-\eta_k^l|\le \min\left\{c\rho^{2+\alpha},c\varepsilon^\frac{1}{4}\rho^2\right\}$. With (\ref{E31}) we get in the same way that $\int_{L_l}\chi_{\phantom{}_{U\cap B^e_\rho(\xi)}}(x,u_\xi^l(x))|\eta^l-\D u_\xi^l(x)|\le \min\left\{c\rho^{2+\alpha},c\varepsilon^\frac{1}{4}\rho^2\right\}.$ Now since $\eta_k^l\to\eta^l$ strongly, we finally get that
\begin{equation}
\liminf_{k\to\infty}(2)\le \min\left\{c\rho^{2+\alpha},c\varepsilon^\frac{1}{4}\rho^2\right\}.
\end{equation}
It remains to estimate the last term (3). It follows as above that
\begin{eqnarray*}
(3) & \le & c\int_{B^e_\rho(\xi)\cap L_l}\left|C_{ij}(x,\overline{u}_k^l(x))-C_{ij}(x,u_\xi^l(x))\right|+c\int_{L_l}\chi_{\phantom{}_{U\cap B^e_\rho(\xi)}}(x,u_\xi^l(x))\left|\D\overline{u}_k^l(x)-\D u_\xi^l(x)\right|. 
\end{eqnarray*}
The first term goes to 0 as usual, and the second term is the same as above, which yields
\begin{equation}
\liminf_{k\to\infty}(3)\le \min\left\{c\rho^{2+\alpha},c\varepsilon^\frac{1}{4}\rho^2\right\}.
\end{equation}
After all we have finally shown that
\begin{equation}
{\cal H}^2\llcorner\left(\graph\overline{u}_k^l\cap B^e_\rho(\xi)\right)(U)={\cal H}^2\llcorner\left(\graph u_\xi^l\cap B^e_\rho(\xi)\right)(U)+\tilde\theta_k(U),
\end{equation}
where $\tilde\theta_k$ is a signed measure such that the $\liminf_{k\to\infty}$ of the total mass is smaller that $\min\left\{c\rho^{2+\alpha},c\varepsilon^\frac{1}{4}\rho^2\right\}$. After passing to a subsequence, $\tilde\theta_k$ converges weakly to some signed measure $\tilde\theta_\xi$ with total mass smaller than $\min\left\{c\rho^{2+\alpha},c\varepsilon^\frac{1}{4}\rho^2\right\}$. Assume that $\tilde\theta_\xi(\partial U)=0$. It follows that $\tilde\theta_k(U)\to\tilde\theta_\xi(U)$, and therefore we get that
\begin{equation}
\lim_{k\to\infty}{\cal H}^2\llcorner\left(\graph\overline{u}_k^l\cap B^e_\rho(\xi)\right)(U)={\cal H}^2\llcorner\left(\graph u_\xi^l\cap B^e_\rho(\xi)\right)(U)+\tilde\theta_\xi(U).
\end{equation}
3.) Since the $\theta_k$'s were signed measures such that the $\liminf$ of the total mass $\le\min\left\{c\rho^{2+\alpha},c\varepsilon\rho^2\right\}$, they converge in the weak sense (after passing to a subsequence) to a signed measure $\overline{\theta}_\xi$ with total mass smaller than $\min\left\{c\rho^{2+\alpha},c\varepsilon\rho^2\right\}$. Assuming $\overline{\theta}_\xi(\partial U)=0$, it follows that $\theta_k(U)\to\overline{\theta}_\xi(U)$.\\
\\
Now by taking limits in (\ref{E19}) we get from 1.), 2.) and 3.) that
$$\mu\llcorner B^e_\rho(\xi)(U)=\sum_{l=1}^M{\cal H}^2\llcorner\left(\graph u_\xi^l\cap B^e_\rho(\xi)\right)(U)+\theta_\xi(U),$$ 
where $\theta_\xi=\overline{\theta}_\xi+\tilde\theta_\xi$ is a signed measure with total mass smaller than $\min\left\{c\rho^{2+\alpha},c\varepsilon^\frac{1}{4}\rho^2\right\}$. Notice that this equation holds for all open $U\subset M$ with $\mu\llcorner B^e_\rho(\xi)(\partial U)={\cal H}^2\llcorner(\graph u_\xi^l\cap B^e_\rho(\xi))(\partial U)=\overline\theta_\xi(\partial U)=\tilde\theta_\xi(\partial U)=0.$ By choosing an appropriate exhaustion this equation holds for arbitrary open sets $U\subset M$ and (\ref{E20}) is shown.\\
\\
Next we claim that $\spt\mu$ is locally given by the union of the graphs of the functions $u_\xi^l$, i.e. for $\rho\le\frac{\rho_0}{256}$ it follows that
\begin{equation}\label{E21}
\spt\mu\cap B_{\rho}^e(\xi)=\bigcup_{l=1}^M\graph u_\xi^l\cap B_{\rho}^e(\xi).
\end{equation}
To prove this let $\rho\le\frac{\rho_0}{128}$ such that \eqref{E20} holds.\\
\\
1.) Let $x\in\spt\mu\cap B_{\frac{\rho}{2}}^e(\xi)$. Proposition \ref{pro:LowerDensity} yields $\mu\llcorner B_\rho^e(\xi)(B_{\frac{\rho}{2}}^e(x))=\mu(B_{\frac{\rho}{2}}^e(x))\ge c\rho^2.$ We get 
$$c\rho^2\le\sum_{l=1}^M{\cal H}^2\left(\graph u_\xi^l\cap B_{\frac{\rho}{2}}^e(x)\right)+c\varepsilon^\frac{1}{4}\rho^2.$$
By choosing $\varepsilon\le\varepsilon_0$ we conclude that $\sum_{l=1}^M{\cal H}^2\left(\graph u_\xi^l\cap B_{\frac{\rho}{2}}^e(x)\right)>0$ and therefore  $x\in\bigcup_{l=1}^M\graph u_\xi^l$.\\
\\
2.) Let $z\in\bigcup_{l=1}^M\graph u_\xi^l\cap B_{\frac{\rho}{2}}^e(\xi)$. Write $z=x+u_\xi^l(x)$ for some $l\in\{1,\ldots,M\}$ and some $x\in L_l$. If $y\in B^e_{\frac{\rho}{4}}(x)\cap L_l$ we claim that $y+u_\xi^l(y)\in B_{\frac{\rho}{2}}^e(z)$, indeed for $\varepsilon\le\varepsilon_0$ we get
$$|z-y-u_\xi^l(y)| \le |x-y|+|u_\xi^l(x)-u_\xi^l(y)| \le \left(1+c\varepsilon^\frac{1}{6}\right)|x-y| \le  \left(1+c\varepsilon^\frac{1}{6}\right)\frac{\rho}{4}\le\frac{\rho}{2}.$$
Therefore 
$${\cal H}^2\llcorner\graph u_\xi^l(B_{\frac{\rho}{2}}^e(z))\ge c \Lzwei(B_{\frac{\rho}{4}}^e(x)\cap L_l) =c\rho^2.$$
As above we obtain $\mu(B_{\frac{\rho}{2}}^e(z))\ge c\rho^2-c\varepsilon^\frac{1}{4}\rho^2>0$ for $\varepsilon\le\varepsilon_0$, and conclude that $z\in\spt \mu$.\\
\\
Now \eqref{E21} implies that the functions $u_\xi^l$ do not depend on the point $\xi$ in the following sense: Let $\eta\in\Sigma\cap B_\frac{\rho_0}{2}^e(\xi_0)$. Then we have for all $\rho\le\frac{\rho_0}{256}$ that 
\begin{equation}\label{E24}
\bigcup_{l=1}^M\graph u_\xi^l\cap\left(B_{\rho}^e(\xi)\cap B_{\rho}^e(\eta)\right)=\bigcup_{l=1}^N\graph u_\eta^l\cap\left(B_{\rho}^e(\xi)\cap B_{\rho}^e(\eta)\right).
\end{equation}
In the next step choose $\rho\le\frac{\rho_0}{256}$ such that $\mu\left(\partial B_\rho^e(\xi)\right)={\cal H}^2_g\llcorner\graph u_\xi^l\left(\partial B_\rho^e(\xi)\right)=0$ for all $l$, and that therefore, from \eqref{E20},
\begin{equation}\label{E26}
\mu\llcorner B_\rho^e(\xi)=\sum_{l=1}^M{\cal H}^2\llcorner\left(\graph u_\xi^l\cap B_\rho^e(\xi)\right)+\theta_\xi.
\end{equation}
Let $z\in\spt\mu\cap B_\rho^e(\xi)=\bigcup_{l=1}^M\graph u_\xi^l\cap B_\rho^e(\xi)$ and let $\sigma>0$ such that $B_\sigma^e(z)\subset B_\rho^e(\xi)$ and such that (due to \eqref{E20} for the point $z$) $\mu\llcorner B_\sigma^e(z)=\sum_{l=1}^N{\cal H}^2\llcorner\left(\graph u_z^l\cap B_\sigma^e(z)\right)+\theta_z$, where the total mass of $\theta_z$ is smaller than $c\sigma^{2+\alpha}$. From (\ref{E24}) it follows that $\theta_\xi\left(B_\sigma^e(z)\right)=\theta_z\left(B_\sigma^e(z)\right),$
hence we get a nice decay for the signed measure $\theta_\xi$, namely
\begin{equation}\label{E25}
\lim_{\sigma\to0}\frac{\theta_\xi\left(B_\sigma^e(z)\right)}{\sigma^2}=0\quad\text{for all }z\in\spt\mu\cap B_\rho^e(\xi)=\bigcup_{l=1}^M\graph u_\xi^l\cap B_\rho^e(\xi).
\end{equation}
Now it follows as before that for all $z\in\spt\mu\cap B_\rho^e(\xi)=\bigcup_{l=1}^M\graph u_\xi^l\cap B_\rho^e(\xi)$
\begin{equation}\label{E22}
\liminf_{\sigma\to0}\frac{\sum_{l=1}^M{\cal H}^2\llcorner(\graph u_\xi^l\cap B_\rho^e(\xi))(B_\sigma^e(z))}{\pi\sigma^2}\ge C >0.
\end{equation}
%To prove this let $z\in\spt\mu\cap B_\rho^e(\xi)=\bigcup_{l=1}^M\graph u_\xi^l\cap B_\rho^e(\xi)$ and let $\sigma>0$ be such that $B_\sigma^e(z)\subset B_\rho^e(\xi)$. Let $z\in\graph u_\xi^l$ for some $l$, then 
%\begin{eqnarray*}
%{\cal H}^2_g\llcorner\left(\graph u_\xi^l\cap B_\rho^e(\xi)\right)(B_\sigma^e(z)) & = & {\cal H}^2_g\left(\graph u_\xi^l\cap B_\sigma^e(z)\right)  =  \int_{L_l} \chi_{\phantom{.}_{B_\sigma^e(z)}}(y+u_\xi^l(y)) \; d\mu_g
%\end{eqnarray*}
%Now let $z=x+u_\xi^l(x)$ with $x\in B^e_\rho(\xi)\cap L_l$. We have that
%$$|z-y-u_\xi^l(y)|\le|x-y|+|u_\xi^l(x)-u_\xi^l(y)|\le\left(1+c\varepsilon^\frac{1}{6}\right)|x-y|,$$
%therefore 
%$$\chi_{\phantom{.}_{B_\sigma^e(z)}}(y+u_\xi^l(y))=1\quad\text{if }|x-y|\le\frac{1}{1+c\varepsilon^\frac{1}{6}}\sigma.$$
%Estimating as before (\ref{E22}) follows ($\varepsilon_0$ maybe smaller).
Now for $z\in\spt\mu\cap B_\rho^e(\xi)=\bigcup_{l=1}^M\graph u_\xi^l\cap B_\rho^e(\xi)$ and $\sigma>0$ such that $B_\sigma^e(z)\subset B_\rho^e(\xi)$ it follows from (\ref{E26}), (\ref{E25}) and (\ref{E22}) that
$$\frac{\mu\llcorner B_\rho^e(\xi)\left(B_\sigma^e(z)\right)}{\sum_{l=1}^M{\cal H}^2\llcorner(\graph u_\xi^l\cap B_\rho^e(\xi))\left(B_\sigma^e(z)\right)}=1+\frac{\theta_\xi\left(B_\sigma^e(z)\right)}{\sum_{l=1}^M{\cal H}^2\llcorner(\graph u_\xi^l\cap B_\rho^e(\xi))\left(B_\sigma^e(z)\right)}.$$
Since the right hand side converges to 1, this shows that $\D_{\left({\sum_{l=1}^M{\cal H}^2\llcorner(\graph u_\xi^l\cap B_\rho^e(\xi))}\right)}\left(\mu\llcorner B_\rho^e(\xi)\right)(z)=1$ for all $z\in\spt\mu\cap B_\rho^e(\xi)=\bigcup_{l=1}^M\graph u_\xi^l\cap B_\rho^e(\xi)$. The Lemma now follows from the Theorem of Radon-Nikodym.
\end{pf}\\
Up to now we have shown that, away from the bad points, the limit measure $\mu$ is locally given by $C^{0,1}$-graphs with small gradient bounded by $c\varepsilon^\frac{1}{6}$. In the next step we will show, using the power decay in Lemma \ref{2ff-absch}, that these graphs are actually $C^{1,\alpha}\cap W^{2,2}$-graphs, and that the $L^2$-norm of their Hessians satisfy a similar power decay.

\begin{pro}\label{pro:C1aReg}
For $\varepsilon\le\varepsilon_0$ there exists a $\rho_0=\rho_0(\xi_0,\varepsilon)>0$ such that
\begin{eqnarray*}
& (i) & u^l_{\xi_0} \in C^{1,\alpha}(L_l\cap B^e_{\rh_0}(\xi_0))\cap W^{2,2} (L_l\cap B^e_{\rh_0}(\xi_0)), \\
& (ii) & \int_{B^e_\sigma(x)\cap L_l}|\D^2 u_{\xi_0}^l|^2\le C\sigma^\alpha\quad\text{for all }x\in B^e_{\rho_0}(\xi_0)\cap L_l\text{ and all }\sigma>0\text{ sufficiently small.}
\end{eqnarray*}
\end{pro}

\begin{pf}
By applying Lemma \ref{mu=graph} to $\xi=\xi_0$, we get that for $\varepsilon\le\varepsilon_0$ there exist $\rh_0=\rh_0(\xi_0,\varepsilon_0)>0$, $2$-dimensional planes $L_l\subset T_{\xi_0}M, l=1,\ldots, M_{\xi_0}$, and functions $u^l_{\xi_0} \in C^{0,1}(L_l\cap B^e_{\rh_0}(\xi_0))$ such that for all $\rh\leq \rh_0$
$$\mu\llcorner B_{\rho}^e(\xi_0)=\sum_{l=1}^{M_{\xi_0}}{\cal H}^2\llcorner\left(\graph u_{\xi_0}^l\cap B_{\rho}^e(\xi_0)\right).$$
Because of the uniform bounds on the area and the Willmore energy of the immersions $f_k$ in the induced metric $g_k$, it follows from Lemma \ref{lem:normalcoordinates} that, for $\rh_0$ maybe smaller, we have $\mu^e_k (B^e_{\rh_0}(\xi_0)) \leq C$ and $\int_{B^e_{\rh_0}(\xi_0)} |H^e_k|^2 d\mu^e_k \leq C.$ It follows that $\mu_k^e\llcorner B^e_{\rh_0}(\xi_0)$ defines an integral, rectifiable 2-varifold with uniformly bounded first variation. By a compactness result for varifolds (see \cite{SiGMT}), there exists an integral, rectifiable 2-varifold $\mu^e$ with weak mean curvature vector $H^e\in L^2(\mu)$, such that (after passing to a subsequence) $\mu_k^e\llcorner B^e_{\rh_0}(\xi_0)\to\mu^e$ weakly in the sense of Radon measures and such that
\begin{equation}\label{4'}
\int_U|H^e|^2\,d\mu^e\le\liminf_{k\to\infty}\int_U|H^e_k|^2\,d\mu_k^e\quad\text{for all open }U\subset B^e_{\rh_0}(\xi_0).
\end{equation}
Repeating the proof of Lemma \ref{mu=graph} by replacing everywhere the Hausdorff measure ${\cal H}^2$ of the manifold with the Euclidean Hausdorff measure ${\cal H}^2_e$, we obtain for all $\rho\leq \rh_0$
$$\mu^e\llcorner B_{\rho}^e(\xi_0)=\sum_{l=1}^{M_{\xi_0}} {\cal H}^2_e\llcorner\left(\graph u_{\xi_0}^l\cap B_{\rho}^e(\xi_0)\right).$$
Since the norm of the mean curvature can be bounded by the norm of the second fundamental form, it follows from Lemma \ref{2ff-absch} and the lower semicontinuity above that for all $\xi \in B^e_{\rh_0}(\xi_0)$ and all $\sigma>0$ such that $B^e_\sigma(\xi)\subset B^e_{\rh_0}(\xi_0)$
\begin{equation*}
\int_{B_\sigma^e(\xi)} |H^e|^2 d\mu^e\leq c\liminf_{k\to\infty} \int_{B_\sigma^e(\xi)} |A^e_k|^2 d\mu^e_k \leq c \sigma^\a.
\end{equation*}
By definition of the weak mean curvature and the graph representation of $\mu^e$ it follows that the functions $u_{\xi_0}^l$ are weak solutions to the weak mean curvature equation
$$\sum_{i,j=1}^2\partial_j\left(\sqrt{\det g_l}\,\,g^{ij}_l\partial_i F_l\right)=\sqrt{\det g_l}\,\,H^e\circ F,$$
where $F_l(x)=x+u_{\xi_0}^l(x)$ and $(g_l)_{ij}=\delta_{ij}+\partial_i u_{\xi_0}^l\cdot\partial_j u_{\xi_0}^l$.   \\
\\
Now first of all it follows from a standard difference quotient argument (see \cite{GT}, Theorem 8.8) that $u_{\xi_0}^l\in W^{2,2}_{loc}(L_l\cap B^e_{\rh_0}(\xi_0))$. By applying the weak mean curvature equation to a suitable test function and using the bounds on $\D u_{\xi_0}^l$ and the power decay of the Willmore energy above one gets for $x\in B^e_{\rh_0}(\xi_0)\cap L_l$ and all $\sigma>0$ sufficiently small that
$$\int_{B^e_\frac{\sigma}{2}(x)\cap L_l}|\D^2 u_{\xi_0}^l|^2\le c\int_{B^e_\sigma(x)\setminus B^e_\frac{\sigma}{2}(x)\cap L_l}|\D^2 u_{\xi_0}^l|^2+c\sigma^\alpha.$$
For details see \cite{Schy}. Now again by "hole-filling" we get $\int_{B^e_\frac{\sigma}{2}(x)\cap L_l}|\D^2 u_{\xi_0}^l|^2\le\theta\int_{B^e_\sigma(x)\cap L_l}|\D^2 u_{\xi_0}^l|^2+c\sigma^\alpha$ for some $\theta\in(0,1)$. Applying Lemma \ref{decay} we obtain (ii). Now it follows from a Lemma of Morrey (see \cite{GT}, Theorem 7.19) that
$$\D u_{\xi_0}^l\in C^{0,\alpha}(B^e_{\rho_0}(\xi_0)\cap L_l),$$
and the Lemma is proved.
\end{pf}\\
Therefore we have up to now shown that our limit measure $\mu$ can be written as $C^{1,\alpha}\cap W^{2,2}$-graphs away from the bad points. Now we will handle the bad points $\B$ and prove a similar power decay as in Lemma \ref{2ff-absch} for balls around the bad points. From this decay it will follow that the set of bad points is actually empty. Since the bad points are discrete and since we want to prove a local decay, we assume that there is only one bad point $\xi_0$, and we will again work in normal coordinates around that point.\\
\\
We will start with a technical but useful Lemma.

\begin{lem}\label{lem:perp}
Consider normal coordinates centered in $\xi_0$ on a neighborhood $U\subset M$. For $x \in U$ let $p \in f_k^{-1}(\{x\})$ be a preimage of $x$ and consider the tangent space $T_pf_k$. We denote with $(T_pf_k)^{\perp_e}$ the orthogonal complement in the normal coordinates, and with $\perp_e$ the projection on $(T_pf_k)^{\perp_e}$. Then for every $\epsilon >0$ there exists a $\rh_0=\rho_0(\xi_0,\varepsilon)>0$, such that for $\rho<\rho_0$ and $k$ sufficiently large
\begin{equation}\label{eq:xxieps}
\frac{|(x-\xi_0)^{\perp_e}|_e}{|x-\xi_0|_e}\le\epsilon\quad\text{for all }x\in(\spt\mu_k\cap B^e_\rho(\xi_0)\setminus B^e_\frac{\rho}{2}(\xi_0))\setminus\mathcal{B}_k,
\end{equation}
where $\mathcal{B}_k\subset \spt\mu_k\cap B^e_{\rho_0}(\xi_0)$ with $\mu_k^e(B^e_\rho(\xi_0)\setminus B^e_\frac{\rho}{2}(\xi_0)\cap\mathcal{B}_k)\le c\epsilon\rho^2$.
\end{lem}

\begin{pf}
By Nash's Embedding Theorem we can assume that $M\subset \Rp$ is isometrically embedded for some $p$. Therefore the sequence $\{f_k\}_{k\in \N}$ can also be seen as a sequence of immersions in $\Rp$. Then Proposition \ref{prop:AreaEstEa1} and the uniform bound on the Willmore energy $W(f_k)$ yield $\int |H_{\Sp^2\hookrightarrow \Rp}|^2 d{\cal H}^2_{\Rp} \leq C.$ By (3.32) in \cite{SiL} there exists a $\rho_0>0$ such that for $\rho<\frac{\rho_0}{4}$ and $k$ sufficiently large 
\begin{equation*}
\frac{|(x-\xi_0)^{\perp_{\Rp}}|_{\Rp}}{|x-\xi_0|_{\Rp}}\le\frac{\varepsilon}{2} \quad \text{for all }x\in(f_k(\Sp^2)\cap B^{\Rp}_{2\rho}(\xi_0)\setminus B^{\Rp}_\frac{\rho}{4}(\xi_0))\setminus\mathcal{B}_k,
\end{equation*}
where $\mathcal{B}_k\subset f_k(\Sp^2)\cap B^{\Rp}_\frac{\rho_0}{2}(\xi_0)$ with ${\cal H}^2_{\Rp}(f_k(\Sp^2)\cap B^{\Rp}_{2\rho}(\xi_0)\setminus B^{\Rp}_{\frac{\rho}{4}}(\xi_0)\cap\mathcal{B}_k)\le c\varepsilon\rho^2$.
Now it's easy to see that 
$$\frac{|(x-\xi_0)^{\perp_e}|_e}{|x-\xi_0|_e} \leq \frac{|(x-\xi_0)^{\perp_{\Rp}}|_{\Rp}}{|x-\xi_0|_{\Rp}}+R(\rh),$$
where $R(\rh)\to 0$ as $\rh\to 0$. Therefore, by choosing $\rh_0$ sufficiently small such that for $\rh<\rh_0$ we have  $R(\rh)<\varepsilon/2$,  $M \cap(B^{e}_{\rho}(\xi_0)\setminus B^{e}_{\frac{\rho}{2}}(\xi_0)) \subset M \cap (B^{\Rp}_{2\rho}(\xi_0)\setminus B^{\Rp}_{\frac{\rho}{4}}(\xi_0))$  and $M\cap B^{\Rp}_\frac{\rho_0}{2}(\xi_0)\subset M\cap B^{e}_{\rho_0}(\xi_0) $, we obtain the result.
%for $\rh<\rho_0$
%$$\frac{|(x-\xi_0)^{\perp_e}|_e}{|x-\xi_0|_e}\le\epsilon\quad\text{for all }x\in(f_k(\Sp^2)\cap B^e_\rho(\xi_0)\setminus B^e_\frac{\rho}{2}(\xi_0))\setminus\mathcal{B}_k,$$
%where $\mathcal{B}_k\subset f_k(\Sp^2)\cap B^e_{\rho_0}(\xi_0)$ with ${\cal H}^2_e(f_k(\Sp^2)\cap B^e_\rho(\xi_0)\setminus B^e_\frac{\rho}{2}(\xi_0)\cap\mathcal{B}_k)\le c\epsilon\rho^2$.
\end{pf}\\
Now remember Definition \ref{df:BadPoints} of the bad points. It follows that there exists a $\rho_0=\rho_0(\xi_0,\varepsilon)>0$ such that for $\rho<\rho_0$ and $k$ sufficiently large
\begin{equation*}
\int_{B_{\frac{3}{2}\rho}(\xi_0)\setminus B_\frac{\rho}{4}(\xi_0)}|A_k|^2\,d\mu_k<\frac{\varepsilon^2}{2}.
\end{equation*}
By choosing $\rho_0$ smaller if necessary it follows from Lemma \ref{lem:normalcoordinates} that
\begin{equation}\label{32}
\int_{B_{\rho}^e(\xi_0)\setminus B_\frac{\rho}{2}^e(\xi_0)}|A^e_k|^2\,d\mu^e_k \leq  \varepsilon^2.
\end{equation}
Moreover we get for $\rho<\rho_0$ and $k$ sufficiently large that
\begin{equation}\label{34}
\spt\mu_k\cap\partial B_{\frac{3}{4}\rho}^e(\xi_0)\neq\emptyset.
\end{equation}
To prove this let $\xi_k\in \spt\mu_k$ such that $\xi_k\to\xi_0$. Thus $\spt\mu_k\cap B_{\frac{3}{4}\rho}^e(\xi_0)\neq\emptyset$ for $k$ sufficiently large. Now suppose that $\spt\mu_k\cap\partial B_{\frac{3}{4}\rho}^e(\xi_0)=\emptyset$. Since $\spt\mu_k$ is connected, we get that $\spt\mu_k\subset B_{\frac{3}{4}\rho}^e(\xi_0)$. It follows that 
$$\diam_h(\spt\mu_k)\le c\diam_e(\spt\mu_k)\le c\rho< c\rho_0,$$
and therefore, by choosing $\rho_0$ smaller if necessary, we get a contradiction to the lower diameter bound given in (\ref{lowdiam}).\\
\\
Let $z\in\spt\mu_k\cap\partial B_{\frac{3}{4}\rho}^e(\xi_0)$. Recalling Lemma \ref{lem:LocAreaEst}, we may apply the Graphical Decomposition Lemma to get that
$$\mu_k^e\llcorner\overline{B_{\frac{\rho}{32}}^e(z)}=\sum_{l=1}^{M_k(z)}\mathcal{H}^2_e\llcorner\Big(\Big(\graph u_k^l\cup\bigcup_j P_j^{k,l}\Big)\cap\overline{B_{\frac{\rho}{32}}^e(z)}\Big),$$
where $\Omega_k^l=(B_{\lambda}^e(\pi_{L_k^l}(z))\cap L_k^l)\setminus\bigcup_m d_{k,m}^l$ with $\lambda>\frac{\rho}{16}$, where $L_k^l$ is a 2-dim. plane, and where the sets $d_{k,m}^l\subset L^l_k$ are pairwise disjoint closed discs. We have the following estimates:
\begin{equation}
M_k(z)\le c=c(M),\phantom{\sum_{l,j}\frac{1}{\rho}}
\end{equation}
\begin{equation}
\sum_{l,m}\diam d_{k,m}^l+\sum_{l,j}\diam P_j^{k,l} \le c\varepsilon^\frac{1}{2}\rho,
\end{equation}
\begin{equation}
\frac{1}{\rho}||u_k^l||_{L^\infty(\Omega_k^l)}+||Du_k^l||_{L^\infty(\Omega_k^l)} \le c\varepsilon^{\frac{1}{6}}.
\end{equation}
\begin{rem}\label{comp1}
Notice that $z\in \spt\mu_k\cap\partial B_{\frac{3}{4}\rho}^e(\xi_0)$ was arbitrary. Cover 
$B_{\left(\frac{3}{4}+\frac{1}{128}\right)\rho}^e(\xi_0)\setminus B_{\left(\frac{3}{4}-\frac{1}{128}\right)\rho}^e(\xi_0)$ by finitely many balls $B_{\frac{\rho}{64}}^e$ with center on $\partial B_{\frac{3}{4}\rho}^e(\xi_0)$ and where the number does not depend on $\rho$, namely
$$B_{\left(\frac{3}{4}+\frac{1}{128}\right)\rho}^e(\xi_0)\setminus B_{\left(\frac{3}{4}-\frac{1}{128}\right)\rho}^e(\xi_0)\subset\bigcup_{i=1}^I B_{\frac{\rho}{64}}^e(y_i),$$
where $y_i\in\partial B_{\frac{3}{4}\rho}^e(\xi_0)$ and $I$ is a universal constant. From this it follows that there exist points $\{z_k^1,\ldots,z_k^{J_k}\}\subset \spt\mu_k\cap\partial B_{\frac{3}{4}\rho}^e(\xi_0)$ with $J_k\le I$, such that
\begin{equation}\label{34'}
\spt\mu_k\cap B_{\left(\frac{3}{4}+\frac{1}{128}\right)\rho}^e(\xi_0)\setminus B_{\left(\frac{3}{4}-\frac{1}{128}\right)\rho}^e(\xi_0)\subset\bigcup_{i=1}^{J_k} B_{\frac{\rho}{32}}^e(z_i^k).
\end{equation}
Now denote by 
\begin{equation}\label{34''}
\left\{\Sigma_k^p\big|\,1\le p\le P_k\right\}
\end{equation}
the images via $f_k$ of the connected components of $f_k^{-1} ( B_{\left(\frac{3}{4}+\frac{1}{128}\right)\rho}^e(\xi_0)\setminus B_{\left(\frac{3}{4}-\frac{1}{128}\right)\rho}^e(\xi_0))$. From the above inclusion, the universal bound on $J_k$, the graphical decomposition from above and the universal bound on $M_k(z_i^k)$ we get that 
\begin{equation}\label{34'''}
P_k\le c,
\end{equation}
where $c$ is a universal constant independent on $k$ and $\rho$.
\end{rem}
In the next step we show that
\begin{equation}\label{35}
\dist\left(\xi_0,L_k^l\right)\le c\varepsilon^\frac{1}{6}\rho\quad\text{for all }l\in\{1,\ldots,M_k(z)\}.
\end{equation}
To prove this notice that Proposition \ref{pro:LowerDensity} and Lemma \ref{lem:normalcoordinates} imply
\begin{equation}\label{36}
\mu^e_k(B_\frac{\rho}{32}^e(z))\ge c\rho^2.
\end{equation}
Moreover notice that
$$(\graph u_k^l\cap B_{\frac{\rho}{32}}^e(z))\setminus\mathcal{B}_k\neq\emptyset,$$
where $\mathcal{B}_k$ was defined in Lemma \ref{lem:perp}. This follows from the graphical decomposition above, the diameter estimates for the sets $P_j^{k,l}$, the area estimate concerning the set $\mathcal{B}_k$ and (\ref{36}).\\
\\
Let $y\in(\graph u_k^l\cap B_{\frac{\rho}{32}}^e(z))\setminus\mathcal{B}_k\subset(\spt\mu_k\cap B_\rho^e(\xi_0)\setminus B_\frac{\rho}{2}^e(\xi_0))\setminus\mathcal{B}_k$. It follows that
$$|\xi_0-\pi_{T_yf_k}(\xi_0)|\le\varepsilon|y-\xi_0|\le\varepsilon\left(|y-z|+|z-\xi_0|\right)\le c\varepsilon\rho.$$
Define the perturbed plane $\tilde L_k^l$ by $\tilde L_k^l=L_k^l+(y-\pi_{L_k^l}(y))$. Thus $\dist(\tilde L_k^l,L_k^l)=|y-\pi_{L_k^l}(y)|\le c\varepsilon^\frac{1}{6}\rho$ (since $y\in\graph u_k^l\cap B_{\frac{\rho}{32}}^e(z)$). Now Pythagoras yields $|y-\pi_{\tilde L_k^l}(\pi_{T_yf_k}(\xi_0))|^2\le|y-\pi_{T_yf_k}(\xi_0)|^2\le|y-\xi_0|^2\le c\rho^2$. Since $T_yf_k$ can be parametrized in terms of $\D u_k^l(y)$ over $\tilde L_k^l$, we get that 
$$|\pi_{T_yf_k}(\xi_0)-\pi_{\tilde L_k^l}(\pi_{T_yf_k}(\xi_0))|\le\|\D u_k^l\|_{L^\infty}|y-\pi_{\tilde L_k^l}(\pi_{T_yf_k}(\xi_0))|\le c\varepsilon^\frac{1}{6}\rho.$$
Therefore by triangle inequality we finally get \eqref{35}.\\
\\
%\begin{eqnarray*}
%dist(\xi,L_k^l) & = & \left|\xi-\pi_{L_k^l}(\xi)\right| \\
%& \le & \left|\xi-\pi_{L_k^l}(\pi_{T_y\Sigma_k}(\xi))\right|\\
%& \le & \left|\xi-\pi_{T_y\Sigma_k}(\xi)\right|+\left|\pi_{T_y\Sigma_k}(\xi)-\pi_{\tilde L_k^l}(\pi_{T_y\Sigma_k}(\xi))\right| \\
%& & +\left|\pi_{\tilde L_k^l}(\pi_{T_y\Sigma_k}(\xi))-\pi_{L_k^l}(\pi_{T_y\Sigma_k}(\xi))\right| \\
%& \le & c\varepsilon^\frac{1}{6}\rho,
%\end{eqnarray*}
Since $\dist(\xi_0,L_k^l)\le c\varepsilon^\frac{1}{6}\rho$, we may assume (after translation) that $\xi_0\in L_k^l$ for all $l\in\{1,\ldots,M_k(z)\}$ and $k$ without changing the estimates for the functions $u_k^l$. Moreover we again have that $L_k^l\to L^l$ with $\xi_0\in L^l.$ Therefore for $k$ sufficiently large we may assume that $L_k^l$ is a fixed 2-dim. plane $L^l$. \\
\\
Now we have that either the point $z$ lies in one of the graphs or can be connected to one of the graphs. Without loss of generality we may assume that this graph corresponds to the function $u_k^1$. Subsequently we will work only with this function $u_k^1$, which is defined on some part of the plane $L_1$ with some discs $d_{k,m}^1$ removed. We will therefore drop the index $1$. Define the set
$$T_k(z)=\left\{\tau\in\left(\frac{\rho}{64},\frac{\rho}{\sqrt{2}\cdot32}\right)\,\Bigg|\,\partial B_\tau^e(\pi_{L}(z))\cap\bigcup_m d_{k,m}=\emptyset\right\}.$$
It follows from the diameter estimates and the selection principle in \cite{SiL} that for $\varepsilon\le\varepsilon_0$ there exists a $\tau\in\left(\frac{\rho}{64},\frac{\rho}{\sqrt{2}\cdot32}\right)$ such that $\tau\in T_k(z)$ for infinitely many $k$.\\
\\
Since $\xi_0\in L$, it follows from the choice of $\tau$ that for $\varepsilon\le\varepsilon_0$
$$\partial B_{\frac{3}{4}\rho}^e(\xi_0)\cap\partial B_\tau^e(\pi_{L}(z))\cap L=\left\{p_{1,k},p_{2,k}\right\},$$
where $p_{1,k}, p_{2,k}\in(B_{\frac{\rho}{\sqrt{2}\cdot32}}^e(\pi_{L}(z))\cap L)\setminus\bigcup_m d_{k,m}$ are distinct points. Define the image points $z_{i,k}\in\graph u_k$ by
$$z_{i,k}=p_{i,k}+u_k(p_{i,k}).$$
Using the $L^\infty$-estimates for $u_k$, we get for $\varepsilon\le\varepsilon_0$ that $\frac{5}{8}\rho<|z_{i,k}-\xi_0|<\frac{7}{8}\rho$, thus $\int_{B_\frac{\rho}{8}^e(z_{i,k})}|A^e_k|^2\,d\mu^e_k\le\varepsilon^2.$ Therefore we can again apply the Graphical Decomposition Lemma to the points $z_{i,k}$. Thus we get that 
$$\mu_k^e\llcorner\overline{B_{\frac{\rho}{32}}^e(z_{i,k})}=\sum_{l=1}^{M_{i,k}(z_{i,k})}{\cal H}^2_e\llcorner\Big(\Big(\graph u_{i,k}^l\cup\bigcup_j P_j^{i,k,l}\Big)\cap\overline{B_{\frac{\rho}{32}}^e(z_{i,k})}\Big),$$
where the usual properties and estimates holds. \\
\\
%\begin{enumerate}
%\item The sets $P_n^{i,k}$ are closed topological discs disjoint from $\graph u_{i,k}$.
%\item $u_{i,k}^l\in C^\infty(\overline{\Omega_{i,k}^l},(L_{i,k}^l)^\perp)$, where $L_{i,k}^l$ is a 2-dim. plane and %$\Omega_{i,k}^l=\left(B_{\lambda_{i,k}^l}(\pi_{L_{i,k}^l}(z_{i,k}))\cap L_{i,k}^l\right)\setminus\bigcup_m d_{i,k,m}^l$, where $\lambda_{i,k}^l>\frac{\rho}{16}$ and where the sets $d_{i,k,m}^l$ are pairwise disjoint closed discs in $L_{i,k}^l$.
%\item The following inequalities hold:
%\begin{eqnarray}
%& & \hspace{-1cm}M_{i,k}(z_{i,k})\le c,\quad\text{where }c<\infty\text{ does not depend on }z_{i,k}, k\text{ and }\rho,\phantom{\sum_n\frac{1}{\rho}} \\
%& & \hspace{-1cm}\sum_m\diam d_{i,k,m}^l+\sum_n\diam P_n^{i,k}\le c\varepsilon^\frac{1}{2}\rho,\phantom{\sum_n\frac{1}{\rho}}\\
%& & \hspace{-1cm}\frac{1}{\rho}\|u_{i,k}^l\|_{L^\infty(\Omega_{i,k}^l)}+\|\D u_{i,k}^l\|_{L^\infty(\Omega_{i,k}^l)}\le c\varepsilon^{\frac{1}{6}}.\phantom{\frac{1}{\rho}\sum_n} 
%\end{eqnarray}
%\end{enumerate}  
Now we have again that the points $z_{i,k}$ either lie in one of the graphs $u_{i,k}^l$ or can be connected to one of them. Without loss of generality let this be the graph corresponding to $u_{i,k}^1$. We will again drop the upper index. Since $z_{i,k}\in\graph u_k$ it follows that $\dist(z_{i,k},L)\le c\varepsilon^\frac{1}{6}\rho$ and that $\graph u_{i,k}$ is connected to $\graph u_k$. Since the $L^\infty$-norms of $u_k$ and $u_{i,k}$ and their derivatives are small, we may assume (after translation and rotation as done before) that $L_{i,k}=L$.\\
\\
By continuing with this procedure we get after a finite number of steps, depending not on $\rho$ and $k$, an open cover of $\partial B_{\frac{3}{4}\rho}^e(\xi_0)\cap L$, which also covers the set
$$A(L)=\left\{x+y\,\Big|\,x\in L, \dist\left(x,\partial B_{\frac{3}{4}\rho}^e(\xi_0)\cap L\right)<\frac{\rho}{\sqrt{2}\cdot64}, y\in L^\perp, |y|<\frac{\rho}{\sqrt{2}\cdot64}\right\}.$$
Now it can happen that after one "walk-around" we do not end up in the same disc of $\spt\mu_k\cap B_{\frac{\rho}{32}}^e(z)$ which contains the point $z$. But then we can proceed in a similar way and do another "walk-around". Now by construction, the "flatness" of the involved graph functions and the diameter bounds for the discs, every "walk-around" corresponds to a part of $\spt\mu_k$ with an area that is bounded from below by $c\rho^2$, where $c$ is a universal constant independent of $k$ and $\rho$. On the other hand we have that $\mu^e_k(B_\rho^e(\xi_0)) \leq c\rho^2$. It follows that after a finite number of "walk-arounds" (which is bounded by a universal constant) we have to get back to the disc of $\spt\mu_k\cap B_{\frac{\rho}{32}}^e(z)$ which contains the point $z$. \\
\\
We summarize the above procedure and the resulting properties in the following remark.
\begin{rem}\label{comp2}
If $\varepsilon\le\varepsilon_0$, then for each component $\Sigma_k^p$ there exist a natural number $k_p$ and a smooth function $u_k^p$ defined on the rectangular set
$$B_k^p=\left[\left(\left(\frac{3}{4}-\frac{1}{\sqrt{2}\cdot64}\right)\rho,\left(\frac{3}{4}+\frac{1}{\sqrt{2}\cdot64}\right)\rho\right)\times[0,2\pi k_p)\right]\setminus\bigcup d_{k,m}^p,$$
where the $d_{k,m}^p$ are closed discs in $\left(\left(\frac{3}{4}-\frac{1}{\sqrt{2}\cdot64}\right)\rho,\left(\frac{3}{4}+\frac{1}{\sqrt{2}\cdot64}\right)\rho\right)\times[0,2\pi k_p)$, such that 
$$\Sigma_k^p=\Big(R_p\left(\graph U_k^p\right)\cup\bigcup_j P_j^{k,p}\Big)\cap B_{\left(\frac{3}{4}+\frac{1}{128}\right)\rho}^e(\xi_0)\setminus B_{\left(\frac{3}{4}-\frac{1}{128}\right)\rho}^e(\xi_0),$$
where $\graph U_k^p=\left\{\left(re^{i\theta},u_k^p(r,\theta)\right)\big|\,(r,\theta)\in B_k^p\right\}$
and $R_p$ denotes a rotation such that $R_p(\mathbb{R}^2)=L_p$, where $L_p$ is the 2-dimensional plane with $\xi_0\in L_p$. Moreover we have
$$\sum_m\diam d_{k,m}^p+\sum_j\diam P_j^{k,p}\le c\varepsilon^\frac{1}{2}\rho, \quad \frac{1}{\rho}\|u_k^p\|_{L^\infty(B_k^p)}+\|\D u_k^p\|_{L^\infty(B_k^p)}\le c\varepsilon^{\frac{1}{6}}.$$
We may assume without loss of generality that the discs $d_{k,m}^p$ are pairwise disjoint, since otherwise we can exchange two intersecting discs by one disc whose diameter is smaller than the sum of the diameters of the intersecting discs.
\end{rem}
Now let $\rho\le\rho_0$ and define the set 
$$C_k(\xi_0)=\left\{\sigma\in\left(\left(\frac{3}{4}-\frac{1}{256}\right)\rho,\left(\frac{3}{4}+\frac{1}{256}\right)\rho\right)\Bigg|\,\partial B_\sigma^e(\xi_0)\cap\bigcup_{p,j} P_j^{k,p}=\emptyset, \int_{\partial B_\sigma^e(\xi_0)}|A^e_k|^2\,ds^e_k\le\frac{512}{\rho}\varepsilon^2\right\}.$$
%It follows that
%$$\Leins(C_k(\xi))\ge\frac{1}{512}\rho,$$
%since by the diameter estimates for the "pimples" we have for $\varepsilon\le\varepsilon_0$ that
%$$\Leins\left(\left\{\sigma\in\left(\left(\frac{3}{4}-\frac{1}{256}\right)\rho,\left(\frac{3}{4}+\frac{1}{256}\right)\rho\right)\Bigg|\,\partial B_\sigma^e(\xi)\cap\bigcup_{p,j} P_j^{k,p}=\emptyset\right\}\right)\ge\frac{1}{256}\rho,$$ 
%and therefore we would get, assuming that $\Leins(C_k(\xi))<\frac{1}{512}\rho$,
%\begin{eqnarray*}
%\varepsilon^2 & \ge & \int_{\Sigma_k\cap B_\rho^e(\xi)\setminus B_\frac{\rho}{2}^e(\xi)}|A_e|^2\,d\mu_e \\
%& \ge & \int_{\left\{\sigma\in\left(\left(\frac{3}{4}-\frac{1}{256}\right)\rho,\left(\frac{3}{4}+\frac{1}{256}\right)\rho\right)\big|\,\partial B_\sigma^e(\xi)\cap\bigcup_{p,j} P_j^{k,p}=\emptyset\right\}\setminus C_k(\xi)}\int_{\Sigma_k\cap\partial B_\sigma^e(\xi)}|A_e|^2\,ds_e\,d\sigma\phantom{\int_{B_\frac{\rho}{2}^e}} \\
%& > & \varepsilon^2.\phantom{\int_{B_\frac{\rho}{2}^e}}
%\end{eqnarray*}
Again it follows from the diameter bounds, a simple Fubini argument and Lemma \ref{selection} that there exists a $\sigma\in\left(\left(\frac{3}{4}-\frac{1}{256}\right)\rho,\left(\frac{3}{4}+\frac{1}{256}\right)\rho\right)$ such that $\sigma\in C_k(\xi_0)$ for infinitely many $k\in\mathbb{N}$. For such a $\sigma$ denote by 
\begin{equation}
\left\{\tilde\Sigma _k^q\,|\,1\le q\le Q_k\right\}
\end{equation}
the images of the components of $f_k^{-1}(B_\sigma^e(\xi_0))$. By Remark \ref{comp1}, we get that $Q_k$ is bounded by a universal constant which is independent of $k$ and $\rho$.

\begin{lem}\label{disc}
Suppose that 
$$\frac{1}{2}\int|A^g_k|^2\,d\mu^g_k\le4\pi-\delta$$
for some $\delta>0$ (which holds in our case by Lemma \ref{lem:EaSpr}).  Then for $\varepsilon\le\varepsilon_0$ each $\tilde{\Sigma}_k^q$ is a topological disc, and moreover $k_p=1$ for all $1\le p\le P_k$.
\end{lem}

\begin{pf}
Fix $k\in\mathbb{N}$. First of all we construct a new immersed surface $\bar\Sigma_k$ such that ($\bar{\mu}_k$ denotes the associated Radon measure)
\begin{eqnarray*}
& (i) & \bar\mu_k\llcorner B_\sigma^e(\xi_0)= \mu^e_k \llcorner B_\sigma^e(\xi_0),\phantom{\Bigg|\int_{B_{\left(\frac{3}{4}+\frac{1}{128}\right)\rho}^e(\xi_0)}\Bigg|} \\
& (ii) & \Bigg|\int_{B_{\left(\frac{3}{4}+\frac{1}{128}\right)\rho}^e(\xi_0)\setminus B_\sigma^e(\xi_0)} K_g\,d\bar{\mu}_k\Bigg|\le c\varepsilon^\frac{1}{3},\quad\text{where }K_g=\text{sectional curvature of }\bar{\Sigma}_k, \\
& (iii) & \int_{M\setminus B_{\left(\frac{3}{4}+\frac{1}{128}\right)\rho}^e(\xi_0)}K_g\,d\bar\mu_k=0.\phantom{\Bigg|\int_{B_{\left(\frac{3}{4}+\frac{1}{128}\right)\rho}^e(\xi_0)}\Bigg|}
\end{eqnarray*}
To define $\bar\Sigma_k$ recall Remark \ref{comp2} and notice that $\sum_{p,m}\diam d_{k,m}^p\le c\varepsilon^\frac{1}{2}\rho.$
Now denote by $M_k$ the number of all discs $d_{k,m}^p$. Because of the diameter estimate it follows for $\varepsilon\le\varepsilon_0$ that there exists an interval $I_k^p\subset\left(\left(\frac{3}{4}-\frac{1}{256}\right)\rho,\left(\frac{3}{4}+\frac{1}{128}\right)\rho\right)$ with $\Leins(I_k^p)\ge\frac{1}{512M_k}\rho$, such that $\left(I_k^p\times[0,2\pi k_p)\right)\cap\bigcup_m d_{k,m}^p=\emptyset.$\\
\\
Let $I_k^p=(a_k^p,b_k^p)$ and $\varphi_p\in C^\infty((0,\infty)\times[0,2\pi k_p))$ with $0\le\varphi_p\le 1$ such that
$$\varphi_p=1\text{ on }(0,a_k^p)\times[0,2\pi k_p),\quad\varphi_p=0\text{ on  }(b_k^p,\infty)\times[0,2\pi k_p), \quad |\D\varphi_p|\le\frac{c}{\rho}\text{ and }|\D^2\varphi_p|\le\frac{c}{\rho^2}. $$
Now define new "components" $\bar\Sigma_k^p$ by
$$\bar\Sigma_k^p=\Big(\Big(R_p\left(\graph\bar U_k^p\right)\cup\bigcup_j P_j^{k,p}\Big)\cap B_{\left(\frac{3}{4}+\frac{1}{128}\right)\rho}^e(\xi_0)\setminus B_{\left(\frac{3}{4}-\frac{1}{128}\right)\rho}^e(\xi_0)\Big)\cup \Big(L_p\setminus B_{\left(\frac{3}{4}+\frac{1}{128}\right)\rho}^e(\xi_0)\Big),$$
where $\graph\bar U_k^p$ is given by
$$\graph\bar U_k^p=\left\{\left(re^{i\theta},\varphi_p(r,\theta)u_k^p(r,\theta)\right)\big|\,(r,\theta)\in B_k^p\right\},$$
and where again $R_p$ denotes a rotation such that $R_p(\mathbb{R}^2)=L_p$. Namely we just "flattened out" the components $\Sigma_k^p$. Observe that by construction and Remark \ref{comp2}, outside of the ball $B_{\left(\frac{3}{4}+\frac{1}{128}\right)\rho}^e(\xi_0)$, $\bar\Sigma_k^p$ is a $k_p$-fold covering of the plane $L_p$.\\
\\
Now define the new surface $\bar\Sigma_k$ by
$$\bar\Sigma_k=\Big(\left(f_k(\Sp^2)\cap B_{\left(\frac{3}{4}+\frac{1}{128}\right)\rho}^e(\xi_0)\right)\setminus\bigcup_p \Sigma_k^p\Big)\cup\bigcup_p\bar\Sigma_k^p.$$
Observe that by construction, $\bar\Sigma_k$ is an immersed surface given by an immersion $F_k:N_k\to M$, where $N_k$ is obtained by gluing ends $E_k^p\cong \R^2\setminus D$ to $f_k^{-1}(B_{\left(\frac{3}{4}-\frac{1}{128}\right)\rho}^e(\xi_0))$ along $f_k^{-1}(\partial\Sigma_k^p\cap\partial B_{\left(\frac{3}{4}-\frac{1}{128}\right)\rho}^e(\xi_0))$, such that outside the ball $B_{\left(\frac{3}{4}+\frac{1}{128}\right)\rho}^e(\xi_0)$, ${F_k}_{|_{E_k^p}}$ is a $k_p$-fold covering of the plane $L_p$.  \\
\\
By definition, $(i)$ follows immediately. Since $\bar\Sigma_k\setminus B_{\left(\frac{3}{4}+\frac{1}{128}\right)\rho}^e(\xi_0)=\bigcup_p L_p\setminus B_{\left(\frac{3}{4}+\frac{1}{128}\right)\rho}^e(\xi_0)$, also (iii) follows directly. To prove property (ii) notice that
\begin{eqnarray*}
\int_{B_{\left(\frac{3}{4}+\frac{1}{128}\right)\rho}^e(\xi_0)\setminus B_\sigma^e(\xi_0)}|K_g|\,d\bar\mu_k & \le & \int_{B_\rho^e(\xi_0)\setminus B_\frac{\rho}{2}^e(\xi_0)}|K_g|\,d\mu^e_k+\sum_p\int_{R_p\left(\graph\bar U_k^p\right)}|K_g|\,d\bar\mu_k.
\end{eqnarray*}
Now the first integral on the right hand side can be estimated by 
$$\int_{B_\rho^e(\xi_0)\setminus B_\frac{\rho}{2}^e(\xi_0)}|K_g|\,d\mu^e_k\le\frac{1}{2}\int_{ B_\rho^e(\xi_0)\setminus B_\frac{\rho}{2}^e(\xi_0)}|A^e_k|^2\,d\mu^e_k\le\varepsilon^2.$$
The second integral can be estimated by
$$\int_{R_p\left(\graph\bar U_k^p\right)}|K_g|\,d\bar\mu_k\le\frac{1}{2}\int_{\graph\bar U_k^p}|A_e|^2\,d\bar\mu_k\le c\int_{B_k^p}|\D^2(\varphi_p u_k^p)|^2.$$
Because of the properties of the functions $u_k^p$ and $\varphi_p$ we have
$$|\D^2(\varphi_p u_k^p)|^2\le c\left(|u_k^p|^2|\D^2\varphi_p|^2+|\D u_k^p|^2|\D\varphi_p|^2+|\varphi_p|^2|\D^2 u_k^p|^2\right)\le c\frac{\varepsilon^\frac{1}{3}}{\rho^2}+|\D^2 u_k^p|^2,$$
and therefore we get
$$ \int_{B_k^p}|\D^2(\varphi_p u_k^p)|^2 \le  c\varepsilon^\frac{1}{3}+c\int_{\graph U_k^p}|A^e_k|^2\,d\mu^e_k \le  c\varepsilon^\frac{1}{3}+c\int_{B_\rho^e(\xi_0)\setminus B_\frac{\rho}{2}^e(\xi_0)}|A^e_k|^2\,d\mu^e_k 
 \le  c\varepsilon^\frac{1}{3}.$$
Thus property (ii) follows by summing over $1\le p\le P_k\le c$.\\
\\
Now denote by $N:\bar\Sigma_k\to\mathbb{S}^2$ the Gau{\ss}-map and notice that $N$ is constant on each end. Therefore the degree of the Gau{\ss}-map $\deg(N)$ is half the Euler characteristic, and it follows from Gau{\ss}-Bonnet that
$$ \deg(N) =  \frac{1}{4\pi}\int K_g\,d\bar\mu_k  =  \frac{1}{4\pi}\int_{B_{\left(\frac{3}{4}+\frac{1}{128}\right)\rho}^e(\xi_0)\setminus B_\sigma^e(\xi_0)}K_g\,d\bar\mu_k+\frac{1}{4\pi}\int_{B_\sigma^e(\xi_0)}K_g\,d\bar\mu_k.$$
Therefore we get that, using (ii) above,
\begin{equation*}\label{eq:degBar}
\Bigg|\int_{B_\sigma^e(\xi_0)}K_g\,d\bar\mu_k-4\pi\deg(N)\Bigg|\le c\varepsilon^\frac{1}{3}.
\end{equation*}
On the other hand it follows from the assumptions and Lemma \ref{lem:normalcoordinates} that
$$ \left|\int_{B_\sigma^e(\xi_0)}K_g\,d\bar\mu_k\right| =  \left|\int_{B_\sigma^e(\xi_0)}K_g\,d\mu^e_k\right| \le  \frac{1}{2}\int_{B_\sigma^e(\xi_0)}|A^e_k|^2\,d\mu^e_k  \le  4\pi-\frac{\delta}{2}$$
by choosing $\rho_0$ smaller if necessary. Since $\deg(N)\in\mathbb{Z}$, it follows for $\varepsilon\le\varepsilon_0$ that
\begin{equation}
\left|\int_{B_\sigma^e(\xi_0)}K_g\,d\mu^e_k\right|=\left|\int_{B_\sigma^e(\xi_0)}K_g\,d\bar\mu_k\right|\le c\varepsilon^\frac{1}{3}.
\end{equation}
Now by the choice of $\sigma$ we have for all $p=1,\ldots,P_k$ that
$$\Sigma^p_k\cap\partial B_\sigma^e(\xi_0)=\gamma_p,$$
where each $\gamma_p$ is a closed, immersed smooth curve and where $P_k$ is bounded by a universal constant. By construction and the choice of $\sigma$ we have that $\gamma_p\cap\bigcup_j P_j^{k,p}=\emptyset$, therefore (see the almost graph representation of $\Sigma_k^p$ above)  $\gamma_p$ is almost a flat circle of radius $\sigma$ which can be parametrized on the interval $[0,2\pi k_p)$. After some computations it follows from the choice of $\sigma$ that (where $\kappa$ denotes the geodesic curvature)
$$ \left|\int_{\gamma_p}\kappa\,ds^e_k-2\pi k_p\right|  \le  c\varepsilon^\frac{1}{6}+c\int_{\gamma_p}|A^e_k|\,ds^e_k  \le  c\varepsilon^\frac{1}{6}+c\sigma^\frac{1}{2}\left(\int_{\partial B_\sigma^e(\xi_0)}|A^e_k|^2\,ds^e_k\right)^\frac{1}{2}  \le  c\varepsilon^\frac{1}{6}+c\left(\frac{\sigma}{\rho}\right)^\frac{1}{2}\varepsilon  \le  c\varepsilon^\frac{1}{6},$$
and therefore it follows from the bound on $P_k$ that
\begin{equation}
\left|\int_{\partial B_\sigma^e(\xi_0)}\kappa\,ds^e_k-2\pi\sum_{p=1}^{P_k}k_p\right|\le c\varepsilon^\frac{1}{6}.
\end{equation}
Now the Euler characteristic of $\tilde\Sigma_k^q$ is given by
$$\chi(\tilde\Sigma_k^q)=2(1-g_q)-b_q,$$
where $b_q$ is the number of boundary components of $\tilde\Sigma_k^q$ and $g_q$ is the genus of the closed surface which arises by gluing $b_q$ topological discs. Especially we have 
$$b_q\ge1\quad\text{and}\quad\sum_{q=1}^{Q_k}b_q=P_k.$$
By summing over $q$ we get that the Euler characteristic of $\bigcup_{q=1}^{Q_k}\tilde \Sigma_k^q$ is
$$\chi_E\left(\bigcup_{q=1}^{Q_k}\tilde \Sigma_k^q\right)=2(Q_k-g)-P_k, \quad \text{where }g=\sum_{q=1}^{Q_k}g_q\ge0.$$
Since $Q_k\le P_k$, we finally get that
$$P_k \ge  2(Q_k-g)-P_k = \frac{1}{2\pi}\int_{B_\sigma^e(\xi_0)}K_g\,d\mu^e_k+\frac{1}{2\pi}\int_{\partial B_\sigma^e(\xi_0)}\kappa\,ds^e_k\ge \sum_{p=1}^{P_k}k_p-c\varepsilon^\frac{1}{6} \ge P_k-c\varepsilon^\frac{1}{6}.$$
Since $2(Q_k-g)-P_k\in\mathbb{N}$, it follows for $\varepsilon\le\varepsilon_0$ that $P_k=2(Q_k-g)-P_k$. Since $Q_k\le P_k$ we get that $Q_k=P_k$ and $g=0$. Thus $g_q=0$ and $b_q=1$ for all $q$. This yields that the Euler characteristic of $\tilde\Sigma_k^q$ is 1 and therefore each $\tilde \Sigma^q_k$ is a topological disc. Moreover the estimate above yields $k_p=1$.
\end{pf}\\
Now define the sets
$$C_k^p=\left\{s\in\left(0,\frac{\rho}{128}\right)\,\Bigg|\,\left(\left(\frac{3}{4}\rho+s\right)\times[0,2\pi)\right)\cap\bigcup_m d_{k,m}^p=\emptyset\right\},$$
$$D_k^p=\left\{s\in C^p_k\,\Bigg|\,\int_{R_p\left(\graph {U_k^p}_{|_{\left(\left(\frac{3}{4}\rho+s\right)\times[0,2\pi)\right)}}\right)}|A^e_k|^2\,ds^e_k\le\frac{512}{\rho}\int_{\Sigma_k^p}|A^e_k|^2\,d\mu^e_k\right\}.$$
By the diameter estimates for the discs $d_{k,m}^p$, again a simple Fubini-type argument and Lemma \ref{selection} there exists a $s\in\left(0,\frac{\rho}{128}\right)$ such that $s\in D^p_k$ for infinitely many $k$. It follows that $u_k^p$ is defined on the line $\left(\frac{3}{4}\rho+s\right)\times[0,2\pi)$. Now it follows from Lemma \ref{disc} that $R_p\Big(\graph {U_k^p}_{|_{\left(\left(\frac{3}{4}\rho+s\right)\times[0,2\pi)\right)}}\Big)$ divides $f_k(\Sp^2)$ into two topological discs $\Sigma_1^{k,p}, \Sigma_2^{k,p}$, one of them, w.l.o.g. $\Sigma_1^{k,p}$, intersecting $B^e_{\frac{3}{4}\rho}(\xi_0)$. From the estimates for the function $u_k^p$ and the choice of $s$ we get that $\Sigma_1^{k,p}\subset B^e_{\left(\frac{3}{4}+\frac{1}{128}\right)\rho}(\xi_0)$, and Lemma \ref{lem:LocAreaEst} yields  $\mu^e_k(\Sigma_1^{k,p})\le c\rho^2.$\\
\\
According to the Lemma \ref{extension},   let $w_k^p\in C^\infty\left(B^e_{\frac{3}{4}\rho+s}(\xi_0)\cap L_p,L_p^\perp\right)$ be an extension of $R_p(U_k^p)$ restricted to $\partial B^e_{\frac{3}{4}\rho+s}(\xi_0)\cap L_p$. In view of the estimates for $u_k^p$, and thus for $w_k^p$, we get that $\graph w_k^p\subset B^e_{\left(\frac{3}{4}+\frac{1}{128}\right)\rho}(\xi_0).$\\
\\
Now we can define the surface $\tilde\Sigma_k$ by
$$\tilde\Sigma_k=\left(f_k(\Sp^2)\setminus\bigcup_p\Sigma_{1}^{k,p}\right)\cup \bigcup_p \graph w_k^p,$$
and we can do exactly the same as in the proof of Lemma \ref{2ff-absch} to get the same power decay as for the good points, but now for balls around the bad points. But by definition the bad points do not allow a decay like this, and therefore we have proved that there are no bad points. \\
\\
Up to now we have shown that the limit measure $\mu$ is locally given by $C^{1,\alpha}\cap W^{2,2}$-graphs. In the next step we show that there exists a $C^{1,\alpha}\cap W^{2,2}$-immersion $f:\Sp^2 \hookrightarrow M$ such that $\mu$ is the Radon measure associated to this immersion $f$. To prove this we will apply a result of Breuning \cite{Breu}, which involves so called generalized $(r,\lambda)$-immersions (for the Definition see Definition \ref{def:rlImm}).\\
\\
For that recall Lemma \ref{final} and Lemma \ref{mu=graph}, namely for every $\xi\in \spt \mu$ there exist a radius $r_\xi>0$ and a natural number $K_\xi\in\N$ such that
\begin{itemize}
\item[(i)] $\mu_k \llcorner B^e_{r_\xi}(\xi)=\sum_{l=1}^{M_\xi} {\cal H}^2 \llcorner\left(\left(\graph u^k_l \cup \bigcup_j P^{k,l}_j\right) \cap B^e_{r_\xi}(\xi)\right)$ for $k\geq K_\xi$, where $u^k_l$ are smooth functions defined on appropriate planes $L_l$ with the usual properties and estimates,
\item[(ii)] $\mu \llcorner B^e_{r_\xi}(\xi)=\sum_{l=1}^{M_\xi} {\cal H}^2 \llcorner\left(\graph u_l \cap B^e_{r_\xi}(\xi)\right)$, where $u_l$ are $C^{1,\alpha}\cap W^{2,2}$-functions defined on $L_l$. 
\end{itemize}
For $\xi \in \spt \mu$ let $\rho_\xi:=\sup \{r_\xi>0 \text{ such that (i) and (ii) holds}\}$. Since $\spt\mu$ is compact, it follows that $\rh:=\inf \{\rh_\xi: \xi\in \spt \mu\}>0$. Notice that (i) and (ii) holds for $\rho$ instead of $r_\xi$. \\
\\
By compactness of $\spt \mu$ there exist $\{\xi_1,\ldots,\xi_I\}\subset \spt \mu$ such that $\spt \mu \subset \bigcup_{i=1}^I B^{e}_{\frac{\rh}{4}}(\xi_i)$. From the Hausdorff distance sense convergence it also follows that $f_k(\Sp^2) \subset \bigcup_{i=1}^I B^{e}_{\frac{\rh}{4}}(\xi_i)$ for $k$ sufficiently large. \\
\\
Now (i) yields $\mu_k \llcorner B^e_{\rh}(\xi_i)=\sum_{l=1}^{M_{\xi_i}} {\cal H}^2 \llcorner\left(\left(\graph u^{k,i}_l \cup \bigcup_j P^{k,l,i}_j\right)\cap B^e_{\rh}(\xi_i)\right)$. By the diameter estimates for the $P^{k,l,i}_j$ and the selection principle \ref{selection} there exists a $\frac{\bar{\rh}}{2}\in (\frac{\rh}{4},\frac{\rh}{2})$ such that $\partial B_{\bar \rh}(\xi_i)\cap \bigcup_{l,j}P^{k,l,i}_j=\emptyset$ for all $i\in \{1,\ldots,I\}$ and infinitely many $k$.\\
\\
Of course we still have that $f_k(\Sp^2) \cup \spt \mu \subset  \bigcup_{i=1}^I B^{e}_{\frac{\bar \rh}{2}}(\xi_i)$, and also the graphical decomposition as in (i) and (ii) still holds in $B^{e}_{\bar \rh}(\xi_i)$.\\
\\
First consider $f_k(\Sp^2)\cap  B^{e}_{\bar \rh}(\xi_1)$. We replace the pimples $\{P^{k,l,1}_j\}_{l,j}$ of $f_k(\Sp^2)\cap  B^{e}_{\bar \rh}(\xi_1)$ with the extension Lemma \ref{extension} as done in the proof of Lemma \ref{2ff-absch} by graphs of functions with small $C^1$-norms defined on the discs $d_{k,m}^{l,1}$. It follows that the sum of the areas of all these graphs is bounded by $c\sum_m(\diam d_{k,m}^{l,1})^2\le c\varepsilon\bar\rho$, which follows from the diameter estimates for the discs. Notice that by the choice of $\bar \rh$, no pimple intersects $\partial B^{e}_{\bar \rh}(\xi_1)$, and we obtain a new $C^{1,1}$-immersion $f_k^1:\Sp^2 \hookrightarrow M$ such that 
\begin{equation}\label{eq:fk1}
f_k(\Sp^2)\setminus B^{e}_{\bar \rh}(\xi_1)=f_k^1(\Sp^2)\setminus B^{e}_{\bar \rh}(\xi_1), \quad f_k^1(\Sp^2)\cap  B^{e}_{\bar \rh}(\xi_1)= \bigcup_{l=1}^{M_1} \graph w_{l}^{k,1}\cap  B^{e}_{\bar \rh}(\xi_1).
\end{equation}
Moreover the above area estimate yields $\mu_k^1(M)\le \mu_k(M)+c\varepsilon\bar\rho$, where $\mu_k^1$ denotes the Radon measure associated to $f_k^1$. Observe that by construction $w_{l}^{k,1}:L_l^1\cap B^{e}_{\bar \rh}(\xi_1)\to (L_l^1)^\perp$ are $C^{1,1}$-functions satisfying $\frac{1}{\bar \rh} \|w_{l}^{k,1}\|_{L^{\infty}}+ \|\D w_{l}^{k,1}\|_{L^{\infty}}\leq c \varepsilon^{\frac{1}{6}}+\delta_k$, where $\delta_k \to 0$. By construction of the limit graphs representing $\mu$ (see the part after Lemma \ref{2ff-absch} and Lemma \ref{mu=graph}) we have that $w_{l}^{k,1} \to u_{l,1}$ uniformly, where $u_{l,1}$ are the graph functions representing $\mu$, namely $\mu \llcorner B^e_{\bar \rh}(\xi_1)=\sum_{l=1}^{M_1} {\cal H}^2 \llcorner\left(\graph u_{l,1} \cap B^e_{\bar \rho}(\xi_1)\right)$. \\
\\
Now consider a point $\xi_j\in\{\xi_1,\ldots,\xi_I\}$ such that $B^{e}_{\frac{\bar \rh}{2}}(\xi_1)\cap  B^{e}_{\frac{\bar \rh}{2}}(\xi_j)\neq \emptyset$, without loss of generality $j=2$. Recall that $\mu_k \llcorner B^e_{\bar{\rh}}(\xi_2)=\sum_{l=1}^{M_2} {\cal H}^2 \llcorner\left(\left(\graph u^{k,2}_l \cup \bigcup_j P^{k,l,2}_j\right) \cap B^e_{\bar \rh}(\xi_2)\right)$, where $u^{k,2}_l$ are smooth functions defined on appropriate planes $L_l^2$.\\
\\
Observe that $f_k^1(\Sp^2)\cap  B^{e}_{\bar \rh}(\xi_1)\cap  B^{e}_{\bar \rh}(\xi_2)= \bigcup_{l=1}^{M_1} \graph w_{l}^{k,1}\cap  B^{e}_{\bar \rh}(\xi_1)\cap  B^{e}_{\bar \rh}(\xi_2)$, and because of the $C^1$-estimates for $w_{l}^{k,1}$ and $u^{k,2}_l$ and the diameter estimate for the pimples, these functions can be written as graphs over the planes $L_{l}^2$ satisfying analogous estimates. We conclude that  $f_k^1(\Sp^2)\cap  B^{e}_{\bar \rh}(\xi_1)\cap  B^{e}_{\bar \rh}(\xi_2)= \bigcup_{l=1}^{M_1} \graph w_{l}^{k,2}\cap  B^{e}_{\bar \rh}(\xi_1)\cap  B^{e}_{\bar \rh}(\xi_2)$, where now the functions $w_{l}^{k,2}$ are defined on the planes $L_l^2$. From \eqref{eq:fk1}, the graphical representation of $f_k(\Sp^2)\cap B^{e}_{\bar \rh}(\xi_2) \setminus B^{e}_{\bar \rh}(\xi_1)$ and the choice of $\bar \rho$, we can replace the pimples inside $B^{e}_{\bar \rh}(\xi_2) \setminus B^{e}_{\bar \rh}(\xi_1)$ with new graphs as done before obtaining a new $C^{1,1}$-immersion $f_k^2:\Sp^2 \hookrightarrow M$ which is the union of graphs (without pimples) in both balls such that the corresponding graph functions converge uniformly to the graph functions representing $\mu$, and such that $\mu_k^2(M)\le \mu_k(M)+2c\varepsilon\bar\rho$, where $\mu_k^2$ denotes the Radon measure associated to $f_k^1$.\\
\\
Repeating the above procedure $I$ times we obtain a $C^{1,1}$-immersion $\tilde f_k:=f_k^I:\Sp^2 \hookrightarrow M$ such that $\mu_k^I(M)\le \mu_k(M)+Ic\varepsilon\bar\rho$, where $\mu_k^I$ denotes the Radon measure associated to $f_k^I$. Because of the uniform area estimate given in Proposition \ref{prop:AreaEstEa1} we have on particular $\mu_k^I(M)\le C$.\\
\\
Now we show that $\tilde f_k$ is actually a generalized $(r,\lambda)$-immersion. Recall that $\spt \mu \subset  \bigcup_{i=1}^I B^{e}_{\bar \rh}(\xi_i)$ is an open cover of $\spt \mu$. By Lebesgue's Lemma there exists the Lebesgue number $\tilde \rho>0$ such that for every $\xi \in \spt \mu$ we have $B^e_{\tilde \rho}(\xi)\subset B^e_{ \bar \rh}(\xi_i)$ for some $i\in \{1,\ldots, I\}$. Now observe that also $\tilde f_k(\Sp^2)$ converges to $\spt \mu$ in the Hausdorff distance sense (which follows from the uniform convergence of the corresponding graphs), thus $B_{\frac{\tilde \rh}{2}} (\tilde f_k(\Sp^2))\subset \bigcup_{i=1}^I B^{e}_{\bar \rh}(\xi_i)$ for $k$ sufficiently large. Let $p \in \Sp^2$ and observe that $B^e_{\frac{\tilde \rh}{2}}(\tilde f_k(p))\subset B^{e}_{\bar \rh}(\xi_i)$ for some $i$. Therefore by construction of $\tilde f_k$ we have 
$$\tilde f_k (\Sp^2)\cap  B^{e}_{\frac{\tilde \rh}{2}}( \tilde f_k(p)) = \bigcup_{l=1}^{M_i} \graph w_{l}^{k,i}\cap  B^{e}_{\frac{\tilde \rh}{2}}( \tilde f_k(p)),$$ 
where $w^{k,i}_{l}:L_l^i\cap B^{e}_{\bar\rho}( \xi_i)\to (L^l_j)^\perp$ are $C^{1,1}$-functions satisfying $\|\D w^{k,i}_{l}\|_{L^{\infty}}\leq c \varepsilon^{\frac{1}{6}}+\delta_k$, where $\delta_k \to 0$. Now recall that by Nash's embedding theorem we can assume that our ambient manifold $M$ is isometrically embedded in some $\Rp$. Denote by $A^k_{p,{L_l^i}}:\Rp \to \Rp$ an Euclidean isometry which maps the origin to $\tilde f_k(p)$ and the subspace $\Rdue\times \{0\}$ onto $\tilde f_k(p)+(L_l^i-\pi_l^i(\tilde f_k(p)))$, where $\pi_l^i$ denotes the orthogonal projection onto $L_l^i$. We get that 
$$ \tilde f_k (\Sp^2)\cap  B^{e}_{\frac{\tilde \rh}{2}}( \tilde f_k(p)) = \bigcup_{l=1}^{M_i} A^k_{p,{L_l^i}}(\graph \tilde w_{l}^{k,i}\cap  B^{e}_{\frac{\tilde \rh}{2}}(0)),$$ 
where $\tilde w_{l}^{k,i}:\Rdue\cap B^{e}_{\bar\rho}( 0)\to (\Rdue)^\perp$ are $C^{1,1}$-functions given by 
$$\tilde w_{l}^{k,i}(x)=\left(A^k_{p,{L_l^i}}\right)^{-1}\left(w_{l}^{k,i}\left(A^k_{p,{L_l^i}}(x)-\Big(\tilde f_k(p)-\pi_l^i(\tilde f_k(p))\Big)\right)-\Big(\tilde f_k(p)-\pi_l^i(\tilde f_k(p))\Big)\right),$$ 
and which satisfy $\|\D \tilde w_{l}^{k,i}\|_{L^{\infty}}\leq c \varepsilon^{\frac{1}{6}}+\delta_k$, where $\delta_k \to 0$. Now denote by $U^k_{\frac{\tilde \rho}{2},p}\subset \Sp^2$ the component of $(\pi\circ (A^k_{p,{L_l^i}})^{-1}\circ \tilde f_k)^{-1}(\Rdue \cap B^{e}_{\frac{\tilde \rh}{2}}(0))$ containing $p$, where $\pi:\Rp\to \Rdue$ is the projection on the first two coordinates. By construction we have $(A^k_{p,{L_l^i}})^{-1}\circ \tilde f_k(U^k_{\frac{\tilde \rho}{2},p})=\graph \tilde w_{l}^{k,i}\cap  B^{e}_{\frac{\tilde \rh}{2}}(0)$ for some $l\in\{1,\ldots,M_i\}$.\\ 
\\
Finally, for given $\lambda<\frac{1}{4}$ we get that for $\varepsilon\le\varepsilon_0$ and $k$ sufficiently large that $\tilde f_k:\Sp^2 \hookrightarrow M$ is a generalized $(\frac{\tilde \rh}{2}, \lambda)$-immersion, namely $\{\tilde f_k\}_{k\in \N}\subset {\cal{F}}^1_C(r,\lambda)$ for $r=\frac{\tilde \rh}{2}$ and $\lambda<\frac{1}{4}$.    \\
\\
By the compactness Theorem \ref{thm:CompGenImm} for generalized $(r,\lambda)$-immersions of Breuning \cite{Breu}, there exist a generalized $(\frac{\tilde \rh}{2}, \lambda)$-function $f:\Sp^2 \hookrightarrow M$ (see Definition \ref{def:rlImm}) and diffeomorphisms $\phi_k:\Sp^2 \to \Sp^2$ such that $\tilde f_k \circ \phi_k\to f$ uniformly. Let us briefly recall Breuning's construction of the limit $f$: Let $q\in\Sp^2$ and $q_k=\phi_k(q)$. By the uniform convergence of $\tilde f_k\circ \phi_k$ we have that for $k$ sufficiently large $B^e_{\frac{\tilde \rh}{2}}(\tilde f_k(q_k))\subset B^{e}_{\bar \rh}(\xi_i)$ for some $i$. By the construction above we know that for each large $k$ 
$$(A^k_{q_k,{L_l^i}})^{-1}\circ \tilde f_k(U^k_{\frac{\tilde \rho}{2},q_k})=\graph \tilde w_l^{k,i}\cap  B^{e}_{\frac{\tilde \rh}{2}}(0)\quad\text{for some }l\in\{1,\ldots,M_i\}.$$
Now Breuning proves that there exist $\lambda$-Lipschitz functions $\tilde u_l^i$ such that 
\begin{equation}\label{breu}
\tilde w_l^{k,i}\to \tilde u_l^i \quad,\quad A^k_{q_k,{L_l^i}}\to A_{q,{L_l^i}}\quad \text{and} \quad (A_{q,{L_l^i}})^{-1}\circ  f(U_{\frac{\tilde \rho}{2},q})=\graph \tilde u_l^i\cap  B^{e}_{\frac{\tilde \rh}{2}}(0).
\end{equation}
On the other hand we know from the representation of the limit measure $\mu$ that
$$\mu \llcorner B^e_{\bar \rho}(\xi_i)=\sum_{l=1}^{M_i} {\cal H}^2 \llcorner\left(\graph u_l^i \cap B^e_{\bar \rho}(\xi_i)\right),$$
where $u_l^i$ are $C^{1,\alpha}\cap W^{2,2}$-functions defined on the planes $L_l^i$. By construction of these limit graphs as carried out before we get that the function given by 
$$A^k_{q_k,{L_l^i}}\left(\tilde w_l^{k,i}\left(\hspace{-0,1cm}\left(A^k_{q_k,{L_l^i}}\right)^{-1}\Big(x+\tilde f_k(q_k)-\pi_l^i(\tilde f_k(q_k))\Big)\hspace{-0,1cm}\right)\right)+\Big(\tilde f_k(q_k)-\pi_l^i(\tilde f_k(q_k))\Big)$$
converges uniformly to $u_l^i$. Since $\tilde f_k(q_k)=\tilde f_k(\phi_k(q))\to f(q)$, it follows from (\ref{breu}) that 
$$\tilde u_l^i(x)=\left(A_{q,{L_l^i}}\right)^{-1}\left(u_l^i\left(A_{q,{L_l^i}}(x)-\Big(f(q)-\pi_l^i(f(q))\Big)\right)-\Big(f(q)-\pi_l^i(f(q))\Big)\right).$$
Therefore the function $\tilde u_l^i$ is actually $C^{1,\alpha}\cap W^{2,2}$ and $A_{q,{L_l^i}}(\graph \tilde u_l^i)=\graph u_l^i$. Thus
$$f(U_{\frac{\tilde \rho}{2},q})=A_{q,{L_l^i}}(\graph \tilde u_l^i\cap  B^{e}_{\frac{\tilde \rh}{2}}(0))=\graph u_l^i\cap  B^{e}_{\frac{\tilde \rh}{2}}(f(q)).$$
We have therefore shown that the generalized $(\frac{\tilde \rh}{2}, \lambda)$-function $f:\Sp^2 \hookrightarrow M$ is actually a $C^{1,\alpha}\cap W^{2,2}$-immersion and that $\mu$ is the Radon measure associated to the immersion $f$.\\
\\
Finally we show that $f$ satisfies the Euler-Lagrange equation and is smooth. First we prove 
\begin{equation}
\label{minimizer}
E(f) = \inf\{E(F)| F \in C^1\cap W^{2,2}(\Sp^2,M) \text{ immersed}\}. 
\end{equation}
A standard approximation argument implies that the right hand side equals the 
infimum $\inf_{[\Sp^2,M]} E(f)$ among smooth immersions. Therefore, \eqref{minimizer}
follows if we prove the lower semicontinuity of the functional, i.e.
\begin{equation}
\label{lowersemi}
E(f) \leq \liminf_{k\to\infty} E(f_k).
\end{equation}
For this we employ results about curvature varifolds due to Hutchinson \cite{Hu1}. 
For convenience of the reader, we recall the main points. For an open 
set $U \subset \R^n$, let $f \in C^1 \cap W^{2,2}_{loc}(\Sigma,U)$ be a 
properly immersed surface with induced metric $g$. For any vector field 
$Y \in C^1_c(\Sigma,\R^n)$ we have the first variation formula 
\begin{equation}
\label{eqvariation}
\int_\Sigma \diver_g Y\,d\mu_g = - \int_{\Sigma} \langle H,Y \rangle\,d\mu_g.
\end{equation}
The projection $P(p):\R^n \to \R^n$ onto the tangent space $T_p f = Df(p) T_p \Sigma$ is given by 
$$
P = g^{\alpha \beta} \langle \partial_\alpha f,\,\cdot\,\rangle \partial_\beta f
\in C^0 \cap W^{1,2}_{loc}(\Sigma,\R^{n \times n}).
$$
We define a vector-valued bilinear form $\overline{B}(p): \R^n \times \R^n \to \R^n$ 
of class $L^2_{loc}$ by the formula 
$$
\overline{B}(e_i,e_j) = 
g^{\alpha \beta} \langle \partial_\alpha f, e_i \rangle (\partial_\beta P) \cdot e_j.
$$ 
Note that $\overline{B}(\partial_\alpha f,\,\cdot\,) = \partial_\alpha P$. Now 
we take  $Y = X \circ G_f$ in \eqref{eqvariation}, where 
$X \in C^1(U \times \R^{n \times n},\R^n)$ has compact support in the 
first variable and $G_f(p) = (f(p),P(p))$ is the Gau{\ss} map. Compute
\begin{eqnarray*}
\diver_g Y & = & g^{\alpha \beta} \langle \partial_\alpha Y, \partial_\beta f \rangle\\
& = & 
g^{\alpha \beta} \langle (D_x X) \circ G_f\, \partial_\alpha f, \partial_\beta f \rangle
+ g^{\alpha \beta} \langle (D_P X) \circ G_f\, \partial_\alpha P, \partial_\beta f \rangle\\ 
& = & {\rm tr}^{Tf} (D_x X) \circ G_f  
+ (\partial_{P^k_j} X^i) \circ G_f\,\overline{B}_{ij}^k.
\end{eqnarray*}
The integral $2$-varifold $V_f$ induced by $f$ on $G_2(U) = U \times G(2,p)$ has 
weight measure $\mu_f = {\cal H}^2 \llcorner \theta_f$, where $\theta_f(x) = \# f^{-1}\{x\}$
is the multiplicity function, and we have 
$$
\int_{G_2(U)} \phi(x,P)\,dV_f(x,P) = \int_U \phi(x,T_x\mu_f)\,d\mu_f(x) \quad \mbox{ for all }
\phi \in C^0_c(G_2(U)). 
$$
Following Section 5.2 in \cite{Hu1}, we show that $V_f$ has generalized curvature
given by 
$$
B(x) = \frac{1}{\theta_f(x)} \sum_{p \in f^{-1}\{x\}} \overline{B}(p) \quad 
\mbox{ for } x \in f(\Sigma).
$$
We put $B = 0$ outside $f(\Sigma)$. To prove the claim we must verify that 
\begin{equation}
\label{weakcurvature}
\int_{G_2(U)} {\rm tr}^P (D_x X)\,dV_f + \int_{G_2(U)} \partial_{P^k_j} X^i \,B_{ij}^k \,dV_f
= - \int_{G_2(U)} \langle B_{ii},X \rangle\,dV_f.
\end{equation}
This will follow from the first variation identity above, recalling that 
$T_x\mu_f$ exists for $\mu_f$-almost every $x \in U$, and hence 
$T_x \mu_f = T_p f$ for $p \in f^{-1}\{x\}$. We compute
\begin {eqnarray*}
\int_{G_2(U)} {\rm tr}^P (D_x X)(x,P)\,dV_f(x,P) %& = &
%\int_U {\rm tr}^{T_x \mu_f} (D_x X)(x,T_x\mu_f)\,\theta_f(x) \,d{\cal H}^2(x)\\
%& = & \int_U \sum_{p \in f^{-1}\{x\}} {\rm tr}^{T_p f} (D_x X)(G_f(p))\,d{\cal H}^2(x)\\
& = & \int_\Sigma {\rm tr}^{T_p f} (D_x X)(G_f(p))\,d\mu_g(p).
\end {eqnarray*} 
Secondly,
\begin {eqnarray*}
\int_{G_2(U)} (\partial_{P^k_j} X^i)(x,P) B_{ij}^k(x)\,dV_f(x,P) %& = & 
%\int_U (\partial_{P^k_j} X^i)(x,T_x \mu_f) \Big(\sum_{p \in f^{-1}\{x\}} \overline{B}_{ij}^k(p)\Big)\,d{\cal H}^2(x)\\
%& = & \int_U \sum_{p \in f^{-1}\{x\}} (\partial_{P^k_j} X^i)(G_f(p)) \overline{B}_{ij}^k(p)\,d{\cal H}^2(x)\\
& = & \int_\Sigma (\partial_{P^k_j} X^i)(G_f(p)) \overline{B}_{ij}^k(p)\,d\mu_g(p). 
\end{eqnarray*}
Similarly,
\begin {eqnarray*}
\int_{G_2(U)} \langle B_{ii}(x),X(x,P) \rangle\,dV_f(x,P) %& = &
%\int_U \Big \langle \sum_{p \in f^{-1}\{x\}} \overline{B}_{ii}(p),X(x,T_x\mu_f) \Big\rangle\,d{\cal H}^2(x)\\
%& = & \int_U \sum_{p \in f^{-1}\{x\}} \big\langle \overline{B}_{ii}(p),X(G_f(p))\big\rangle\,d{\cal H}^2(x)\\
& = &  \int_\Sigma \big\langle \overline{B}_{ii}(p),X(G_f(p)) \big\rangle\,d\mu_g(p). 
\end{eqnarray*}
To calculate $\overline{B}$, we first observe that $\overline{B}(N,\,\cdot\,) = 0$ if 
$N$ is normal along $f$. We further calculate 
\begin{eqnarray*}
\overline{B}(\partial_\alpha f,\partial_\beta f) 
& = & \partial_\alpha P \cdot \partial_\beta f
= \partial^2_{\alpha \beta} f - P \partial^2_{\alpha \beta} f
= A_{\alpha \beta}\\
\overline{B}(\partial_\alpha f,N) & = & 
\partial_\alpha P \cdot N = - P \partial_\alpha N 
= - g^{\beta \gamma}  \langle \partial_\alpha N, \partial_\beta f \rangle \partial_\gamma f
= g^{\beta \gamma} \langle N,A_{\alpha \beta} \rangle \partial_\gamma f.
\end{eqnarray*}
In particular $\overline{B}_{ii} = H$ which completes the proof of \eqref{weakcurvature}. 
We will also need that $B \in L^1_{loc}(\mu_f)$ is uniquely determined by \eqref{weakcurvature}, 
see \cite{Hu1}, Proposition 5.2.2.\\
\\
Next consider a sequence of varifolds $V_k \to V$ weakly in $G_2(U)$, 
and functions $\psi_k \in L^2(V_k,\R^m)$ with 
$$
C_0:= \lim_{k \to \infty} \|\psi_k\|_{L^2(V_k)} < \infty. 
$$
Define the linear functionals $\Lambda_k:C^0_c(G_2(U),\R^m) \to \R$,
$\Lambda_k(\phi) = \int_{G_2(U)} \langle \phi, \psi_k \rangle\,dV_k$.
Clearly 
$$
|\Lambda_k(\phi)| \leq 
\|\psi_k\|_{L^2(V_k)} V_k({\rm spt\,} \phi)^{\frac{1}{2}} \|\phi\|_{C^0(U)}.
$$
By the Banach-Alaoglu theorem, we have $\Lambda_k \to \Lambda$ in
$C^0_c(G_2(U))'$ for a subsequence, and we get
$$
|\Lambda(\phi)| \leq C_0\, V({\rm spt\,}\phi)^{\frac{1}{2}} \|\phi\|_{C^0(U)} \quad 
\mbox{ for } \phi \in C^0_c(G_2(U),\R^m).
$$
By the theorem of Riesz, the functional $\Lambda$ has a representation 
$$
\Lambda(\phi) = \int_{G_2(U)} \langle \phi, \nu \rangle \,d|\Lambda|,
$$
where $|\Lambda|$ is the variation measure and $\nu:G_2(U) \to \R^m$ 
is Borel measurable with $|\nu| = 1$ almost everywhere with respect to 
$|\Lambda|$. But $|\Lambda|$ is absolutely continuous with respect to $V$, 
hence we have $|\Lambda| = V \llcorner \theta$ for some function
$\theta \in L^1_{loc}(V,\R^{+}_0)$. Put $\psi = \theta \nu \in L^1_{loc}(V,\R^m)$ 
to obtain 
$$
\int_{G_2(U)} \langle \phi,\psi \rangle\,dV = 
\lim_{k \to \infty} \int_{G_2(U)} \langle \phi,\psi_k \rangle\,dV_k
\quad \mbox{ for all } \phi \in C^0_c(G_2(U),\R^m).
$$
Now for any $\phi \in C^0_c(G_2(U),\R^m)$ we can estimate
$$
\Lambda(\phi) = \int_{G_2(U)} \langle \phi,\psi \rangle \,dV 
\leq C_0 \lim_{k \to \infty} \|\phi\|_{L^2(V_k)} = C_0 \|\phi\|_{L^2(V)}.
$$
Thus $\Lambda$ extends continuously to $L^2(V,\R^m)$ and hence $\psi \in L^2(V,\R^m)$.
Moreover, for any $\eta \in C^0_c(U,\R^{+}_0)$ we get by Cauchy-Schwarz
$$
\int_U \eta |\psi|^2\,dV \leq \Big(\int_U \eta |\psi|^2\,dV \Big)^{\frac{1}{2}}
\liminf_{k \to \infty} \Big(\int_U \eta |\psi_k|^2 \,dV_k\Big)^{\frac{1}{2}}.
$$
Canceling we obtain
$$
\int_U \eta |\psi|^2\,dV \leq \liminf_{k \to \infty} \int_U \eta |\psi_k|^2 \,dV_k.
$$
We now return to the setting of immersed surfaces. Let $f_k \in C^1 \cap W^{2,2}_{loc}(\Sigma_k,U)$,
$f \in C^1 \cap W^{2,2}_{loc}(\Sigma,U)$ be properly immersed. Assume that 
$$
\|A_{f_k}\|_{L^2(\Sigma_k)} \leq C_0,\,\quad \mbox{ and } \quad 
V_{f_k} \to V_f \quad \mbox{ as varifolds in } U.
$$
Let us fix a cutoff function $\eta \in C^0_c(U,\R^{+}_0)$. From the above 
we see $|\overline{B}|^2 = 2 |A|^2$ and 
\begin{eqnarray*}
\int_U \eta(x) |B(x)|^2\,d\mu_f & = & 
\int_U \eta(x) \theta_f(x)^{-2} \Big|\sum_{p \in f^{-1}\{x\}} \overline{B}(p) \Big|^2\,d\mu_f(x)\\
& \leq & \int_U \eta(x) \sum_{p \in f^{-1}\{x\}} |\overline{B}(p)|^2\,d{\cal H}^2(x)\\
& = & 2 \int_{\Sigma} \eta \circ f  |A|^2 \,d\mu_g.
\end{eqnarray*}
In order to have equality for $f$ in this argument, we make the technical assumption 
that $f$ is injective. It now follows that $B_{f_k}$ is bounded in
$L^2(V_k)$, and $V_k \llcorner B_{f_k}$ converges to $V \llcorner B$ 
as varifolds, for some $B \in L^2(V)$. Taking limits in \eqref{weakcurvature} 
shows that $V$ has generalized second fundamental form equal to $B$,
hence we have $B = B_{f}$ by uniqueness. We conclude 
\begin{eqnarray*}
\int_\Sigma \eta \circ f |A_f|^2\,d\mu_g & = & 
\frac{1}{2} \int_{G_2(U)} \eta\, |B|^2\,dV_f\\
& \leq & \frac{1}{2} \liminf_{k \to \infty} \int_{G_2(U)} \eta\, |B_{f_k}|^2\,dV_{f_k}\\
& \leq & \liminf_{k \to \infty} \int_{G_2(U)} \eta \circ f_k |A_{f_k}|^2\,d\mu_{g_k}.
\end{eqnarray*}
This proves the local lower semicontinuity of the functional $E(f)$. Finally, 
assume that $f:\Sigma \to M \cap U$ where $M \subset \R^n$ is a $C^2$ 
submanifold. The second fundamental forms in $M$ and in $\R^n$ differ 
only by a first order term, more precisely 
$$
\int_\Sigma |A|^2\,d\mu_g = 
\int_\Sigma |A^{\R^n}|^2\,d\mu_g - \int_\Sigma |A^M_{Tf \times Tf}|^2\,d\mu_g.
$$
Here by $A^{\R^n}$ we mean the second fundamental form in $\R^n$, while $A$ 
now refers to the second fundamental form in $M$. Extending the second 
fundamental form $A^M$ of $M \subset \R^n$ to $TM^\perp$ by zero, we may write 
\begin{eqnarray*}  
\int_\Sigma \eta \circ f |A^M_{Tf \times Tf}|^2\,d\mu_g & = & 
\int_U \eta(x) \sum_{p \in f^{-1}\{x\}} |A^M(x)(P(p)e_i,P(p)e_j|^2\,d{\cal H}^2(x)\\
& = & \int_{G_2(U)} \eta(x) |A^M(x)(Pe_i,P e_j)|^2 \,dV_f(x,P).
\end{eqnarray*}
The last expression is continuous under the convergence $V_{f_k} \to V_f$. 
Therefore the $L^2$ integral of the second fundamental form in $M$ is also 
lower semicontinuous.\\
\\
To prove the lower semicontinuity of the full functional $E(f)$, we cover 
the image of the limit surface by neighborhoods on which we have a local  
graph description (with pimples for the $f_k$). Then we choose a 
subordinate partition of unity and apply the above lower semicontinuity 
statement to each of the graphs. Summing up yields the desired result. \\
\\
Now by construction and lower semicontinuity it follows that the limit immersion $f$ minimizes $E$ among $C^{1}\cap W^{2,2}$-immersions, in particular it satisfies the Euler-Lagrange equation.\\
\\
To compute the Euler-Lagrange equation, let $f \in W^{2,2} \cap C^{1,\alpha}(U,\R^3)$, $f(x) = (x,u(x))$, be a graph given in local 
coordinates in $M$. The functional $E(f)$ is then given by 
$$
E(f) = \frac{1}{2} \int_U \sqrt{\det g}\, g^{\alpha \gamma} g^{\beta \lambda} 
h (P^\perp D_\alpha \partial_\beta f,D_\gamma \partial_\lambda f),
$$
where $h = h_{ij}$ is the Riemannian metric on $M$, and 
\begin{eqnarray*} 
g_{\alpha \beta} &  = & (h \circ f)(\partial_\alpha f,\partial_\beta f),\\
P^\perp & = & {\rm Id} - g^{\alpha \beta} (h \circ f)(\partial_\alpha f,\, \cdot\,) \partial_\beta f,\\ 
D_\alpha \partial_\beta f & = & (0,\partial^2_{\alpha \beta} u) + 
(\Gamma \circ f)(\partial_\alpha f,\partial_\beta f).
\end{eqnarray*}
Here $\Gamma = \Gamma_{ij}^k$ are the Christoffel symbols of $M$. The functional 
thus has the general form 
$$
E(f) = \int_U \Big(A^{\alpha \beta \gamma \lambda}(x,u,Du) \partial^2_{\alpha \beta}u \,\partial^2_{\gamma \lambda} u 
+ B^{\alpha \beta}(x,u,Du) \partial^2_{\alpha \beta} u
+ C(x,u,Du)\Big),
$$
where $A,B,C$ are smooth functions, and specifically for $e_3 = (0,1) \in \R^3$ 
$$
A^{\alpha \beta \gamma \lambda}(x,u,Du) 
= \frac{1}{2} \sqrt{\det g}\, g^{\alpha \gamma} g^{\beta \lambda}\, 
h\big(P^\perp e_3,e_3)\, \partial^2_{\alpha \beta} u\, \partial^2_{\gamma \lambda} u.
$$
We see that a bound for $Du$ implies an ellipticity condition
$$
A^{\alpha \beta \gamma \lambda} \xi_{\alpha \beta} \xi_{\gamma \lambda} 
= \frac{1}{2} \sqrt{\det g}\, \|P^\perp e_3\|_h^2\, \|\xi\|_g^2 \geq \lambda |\xi|^2 > 0.
$$
It is now straightforward to check that the Euler-Lagrange equation satisfies all 
the conditions of Lemma $3.2$ in \cite{SiL}, provided that $Du$ is bounded. 
Hence we get that $u$ belongs locally to $W^{3,2} \cap C^{2,\alpha}$ for some
$\alpha > 0$, and that the $L^2$ integral of $D^3 u$ satisfies a power decay.
As in \cite{SiL} we can refer to \cite{MCB} to conclude that $u$ is in fact smooth. Therefore Theorem \ref{thm:ExEK} is proved.

\section{Proof of Theorem \ref{thm:ExW1}}
The proof of Theorem \ref{thm:ExW1}, namely the problem of minimizing the functional 
$$W_1(f) = \int_{\Sp^2} \Big(\frac{1}{4}|H|^2 +1\Big)\,d\mu_g$$
in the class of immersions $f:\Sp^2 \hookrightarrow M$, where $M$ is a closed, three-dimensional Riemannian manifold with sectional curvature $K^M \leq 2$ and moreover $R^M(\overline{x}) > 6$ for some point $\overline{x} \in M$, is very similar to the proof of Theorem \ref{thm:ExEK}. Here we summarize the different steps of the proof and point out the differences to the proof of Theorem \ref{thm:ExEK}.\\
\\
Again we use the concept of minimizing sequences. Therefore let $f_k:\Sp^2\hookrightarrow M$ be a minimizing sequence of immersed closed surfaces for the functional $W_1$ and denote by $\mu_k$ the Radon measure on $M$ associated to $f_k$. Obviously we have that $\mu_k(M)\le W_1(f_k)\le C$ uniformly in $k$. Therefore there exists a Radon measure $\mu$ on $M$ such that, up to subsequences, 
\begin{equation}
\mu_k \to \mu \quad \text{weakly as Radon measures},
\end{equation}
and as before the monotonicity formula Lemma \ref{lem:Link} yields
\begin{equation}
\spt \mu_k \to \spt \mu \quad \text{in the Hausdorff distance sense}.
\end{equation}
Observe that, since $R^M(\overline{x}) > 6$ for some point $\overline{x} \in M$, it follows similar to Lemma \ref{lem:EaSpr} that 
\begin{equation}\label{W1<4pi}
\inf_{f \in [\Sp^2,M]} W_1(f) < 4\pi.
\end{equation}
Using Lemma \ref{lem:normalcoordinates}, it follows that Proposition \ref{prop:LBdiamEa} also holds for $E$ replaced by $W_1$, which yields that we again have a lower diameter bound, namely
\begin{equation}\label{lowdiamW_1}
\diam_h (\spt \mu)\geq \liminf_k (\diam_h \spt \mu_k) >0.
\end{equation}
The next Lemma states an important upper bound for the functional $E$ in terms of the functional $W_1$ and is a direct consequence of equation (\ref{eq:EQF}).
\begin{lem}\label{lem:AboundedWe}
Let $M$ be a compact Riemannian manifold with sectional curvature $K^M\leq 2$, and let $f:\Sp^2 \hookrightarrow M$ be a smooth immersion. It follows that
$$E(f)\leq 2W_1(f)-4\pi$$
\end{lem}
It follows that $\limsup_{k\to\infty}E(f_k)<4\pi$. Moreover it follows from this uniform upper bound that we can define the bad points with respect to $\varepsilon>0$ as in Definition \ref{df:BadPoints}, and that also Remark \ref{E7}, Lemma \ref{lem:EstGraDec}, the Graphical Decomposition Lemma \ref{final} and the lower density bound in Proposition \ref{pro:LowerDensity} hold in exactly the same way. \\
\\
Now observe that the proof of the power decay of the $L^2$-norm of the second fundamental form in Lemma \ref{2ff-absch} carries over analogously up to equation (\ref{eq:ProvEstAe}) (for the following notation see the proof of Lemma \ref{final}). Now, since $f_k$ is a minimizing sequence for the functional $W_1$, we have that
\begin{equation}\label{eq:WeSStilde}
W_1(\tilde f_k)\ge W_1(f_k)-\varepsilon_k,\quad\text{where }\varepsilon_k\to0.
\end{equation}
Equation (\ref{eq:EQF}) yields (using that the sectional curvature is bounded by compactness of the manifold $M$)
\begin{equation}\label{eq:S-StildeWe}
\sum_{l=1}^{M_k}\int_{\graph w_k^l}|A|^2\,d{\cal H}^2 + c\sum_{l=1}^{M_k} {\cal H}^2(\graph w_k^l)\geq \int_{B^e_{\frac{\rh}{16}}(\xi)} |A_k|^2\,d\mu_k-c \mu_k(B^e_{\frac{\rh}{16}}(\xi))-\varepsilon_k.
\end{equation}
Using that ${\cal H}^2(\graph w_k^l)\le c\rho^2$ by the estimates for $w_k^l$, it follows from Lemma \ref{lem:LocAreaEst} and Lemma \ref{lem:normalcoordinates} that (\ref{SS-Stilde}) holds in this setting. The rest of the proof is again the same as before. This shows that also Lemma \ref{final} holds.\\
\\
Now we can construct the limit graph functions as done before after the proof of Lemma \ref{2ff-absch}, and show in the same way as before that the limit measure $\mu$ is locally (around the good points) given by the sum of the 2-dimensional Hausdorff measure restricted to these limit graphs, namely that Lemma \ref{mu=graph} holds. Observe that also Proposition \ref{pro:C1aReg} holds, thus the limit measure is given by $C^{1,\alpha}\cap W^{2,2}$-graphs away from the bad points. \\
\\
To exclude the bad points, we can do the same as before. Observe that the crucial Lemma \ref{disc} also holds, because by Lemma \ref{lem:AboundedWe} and (\ref{W1<4pi}) the assumption
$$\frac{1}{2}\int|A^g_k|^2\,d\mu^g_k\le4\pi-\delta$$
for some $\delta>0$ are satisfied. Thus $\mu$ is locally given by $C^{1,\alpha}\cap W^{2,2}$-graphs. \\
\\
As before, using \cite{Breu}, it follows that there exists a $C^{1,\alpha}\cap W^{2,2}$-immersion $f:\Sp^2 \hookrightarrow M$ such that $\mu$ is the Radon measure associated to this immersion $f$. To conclude that $f$ is actually smooth, observe that by construction and lower semicontinuity $f$ satisfies the Euler-Lagrange equation for the functional $W_1$. By equation (\ref{eq:EQF}) the functionals $E$ and $W_1$ differ only by a topological constant and a multiple of $K^M_f$, which is a smooth function of $u$, $\D u$ in graph coordinates. Therefore the conditions of Lemma $3.2$ in \cite{SiL} are again satisfied. Hence we get that $u$ belongs locally to $W^{3,2} \cap C^{2,\alpha}$ for some
$\alpha > 0$, and that the $L^2$ integral of $D^3 u$ satisfies a power decay.
As in \cite{SiL} we can refer to \cite{MCB} to conclude that $u$ is in fact smooth. Therefore also Theorem \ref{thm:ExW1} is proved.

\section{Appendix}

\subsection{Some useful Lemmas}
In this subsection we state some useful results we need for proving regularity. Lemma \ref{extension} is an extension result adapted to the cut-and-paste procedure we use and is proved in \cite{Schy}.
\begin{lem}\label{extension}
Let $L$ be a 2-dimensional plane in $\Rn$, $x_0\in L$ and $u\in\C^\infty\left(U,L^\perp\right)$, where $U\subset L$ is an open neighborhood of $L\cap\partial B_\rho(x_0)$. Moreover let $|\D u|\le c$ in $U$. Then there exists a function $w\in\C^\infty(\overline{B_\rho(x_0)},L^\perp)$ with the following properties:
\begin{eqnarray*}
& 1.) & w=u\quad\text{and}\quad\frac{\partial w}{\partial\nu}=\frac{\partial u}{\partial\nu}\quad\text{on }\partial B_\rho(x_0),\phantom{\int_{B_\rho}} \\
& 2.) & \frac{1}{\rho}||w||_{L^\infty(B_\rho(x_0))}\le c(n)\left(\frac{1}{\rho}||u||_{L^\infty(\partial B_\rho(x_0))}+||\D u||_{L^\infty(\partial B_\rho(x_0))}\right),\phantom{\int_{B_\rho}} \\
& 3.) & ||\D w||_{L^\infty(B_\rho(x_0))}\le c(n)||\D u||_{L^\infty(\partial B_\rho(x_0))},\phantom{\int_{B_\rho}} \\
& 4.) & \int_{B_\rho(x_0)}|\D^2w(x)|^2\,\de x\le c(n)\rho\int_{\graph u_{|_{\partial B_\rho(x_0)}}}|A(x)|^2\, d{\cal H}^1,
\end{eqnarray*}
where $d{\cal H}^1$ is the 1-dimensional Euclidean Hausdorff measure.
\end{lem}

The second lemma is a useful selection principle proved in \cite{SiL}.

\begin{lem}\label{selection}
Let $\delta>0, I\subset\R$ a bounded interval and $A_k\subset I, k\in\N,$ measurable sets with $\Leins(A_k)\ge\delta$ for all $k$. Then there exists a set $A\subset I$ with $\Leins(A)\ge\delta$, such that each point $x\in A$ lies in $A_k$ for infinitely many $k$.
\end{lem}
%\begin{pf}
%Define the set $A\subset I$ by
%$$A=\bigcap_{l=1}^\infty\bigcup_{k\ge l}A_k.$$
%Then if $x\in A$ it follows that $x\in\bigcup_{k\ge l}A_k$ for all $l\in\N$ and therefore $x\in A_k$ for infinitely many $k$. %We also have that
%$$\Leins(A)=\lim_{l\to\infty}\Leins\left(\bigcup_{k\ge l}A_k\right)\ge\delta.$$
%Therefore the Lemma is proven.
%\end{pf}

The third lemma is a decay result we need to get the power decay for the $L^2$-norm
of the second fundamental form in Lemma \ref{2ff-absch}.

\begin{lem}\label{decay}
Let $g:(0,b)\to[0, +\infty)$ be a bounded function such that 
$$g\left(x\right)\le\gamma g(2x)+Cx^\alpha\quad\text{for all }x\in\left(0,\frac{b}{2}\right),$$
where $\alpha>0$, $\gamma\in(0,1)$, and $C\ge0$ is a constant. Then there exists a $\beta\in(0,1)$ and a constant $C=C\left(\gamma,\alpha,b,||g||_{L^\infty(0,b)}\right)$ such that
$$g(x)\le Cx^\beta\quad\text{for all }x\in\left(0,b\right).$$
\end{lem}
The last statement is a generalized Poincar\'e inequality proved in \cite{SiL}.
\begin{lem}\label{Poincare}
Let $\mu>0$, $\delta\in\left(0,\frac{\mu}{2}\right)$ and $\Omega=B^{\Rdue}_\mu(0)\backslash E$, where $E\subset \Rdue$ is measurable with $\Leins(p_1(E))\le\frac{\mu}{2}$ and $\Leins(p_2(E))\le\delta$, where $p_1$ is the projection onto the $x$-axis and $p_2$ is the projection onto the $y$-axis. Then for any $f\in C^1(\Omega)$ there exists a point $(x_0,y_0)\in\Omega$ such that 
$$\int_\Omega\left|f-f(x_0,y_0)\right|^2\le C\mu^2\int_\Omega\left|\D f\right|^2+C\delta\mu\sup_\Omega|f|^2$$
where $C$ is an absolute constant.
\end{lem}

\subsection{Definitions and properties of generalized $(r,\lambda)$-immersions}

Here we recall the definitions and properties of generalized $(r,\lambda)$-immersions $f:\Sp^2\hookrightarrow M\subset \Rp$ appearing in \cite{Breu}.\\
We call a mapping $A:\Rp\to\Rp$ an Euclidean isometry, if there is a rotation $R\in SO(p)$ and a translation $T\in\Rp$, such that $A(x)=Rx+T$ for all $x\in\Rp$.\\
For a given point $q\in\Sp^2$ and a given 2-plane $E\in G(p,2)$ let $A_{q,E}:\Rp\to\Rp$ be an Euclidean isometry which maps the origin to $f(q)$ and the subspace $\Rdue\times\{0\}\subset\Rp$ onto $f(q)+E$.\\
Let $U^E_{r,q}\subset\Sp^2$ be the $q$-component of the set $(\pi\circ A^{-1}_{q,E}\circ f)^{-1}(B_r)$, where $\pi:\Rp\to\Rdue$ is the projection on the first two coordinates.

\begin{df}\label{def:rlImm}
An immersion $f:\Sp^2\hookrightarrow M\subset\Rp$ is called a generalized $(r,\lambda)$-immersion, if for each point $q\in \Sp^2$ there is an $E=E(q)\in G(p,2)$, such that $A^{-1}_{q,E}\circ f(U^E_{r,q})$ is the graph of a differentiable function $u:B_r\to(\Rdue)^\perp$ with $\|\D u\|_{C^0(B_r)}\le\lambda$.\\
The set of generalized $(r,\lambda)$-immersions is denoted by ${\cal{F}}^1(r,\lambda)$. Moreover let ${\cal{F}}^1_V(r,\lambda)$ be the set of all immersions $f\in {\cal{F}}^1(r,\lambda)$ such that $\mu_g(\Sp^2)\le V$, where $\mu_g$ is the induced area measure. \\
A continuous function $f:\Sp^2\hookrightarrow M\subset\Rp$ is called a $(r,\lambda)$-function, if for each point $q\in \Sp^2$ there is an $E=E(q)\in G(p,2)$, such that $A^{-1}_{q,E}\circ f(U^E_{r,q})$ is the graph of a Lipschitz function $u:B_r\to(\Rdue)^\perp$ with with Lipschitz constant $\lambda$. The set of $(r,\lambda)$-functions is denoted by ${\cal{F}}^0(r,\lambda)$.
\end{df}

Now we recall the Compactness Theorem in \cite{Breu}, Theorem 0.5.

\begin{thm}\label{thm:CompGenImm}
Let $\lambda\le\frac{1}{4}$. Then ${\cal{F}}^1_V(r,\lambda)$ is relatively compact in ${\cal{F}}^0(r,\lambda)$ in the following sense: Let $f_k:\Sp^2\hookrightarrow M\subset\Rp$ be a sequence in ${\cal{F}}^1_V(r,\lambda)$. Then, after passing to a subsequence, there exists a function $f\in {\cal{F}}^0(r,\lambda)$ and a sequence of diffeomorphisms $\phi_k:\Sp^2\to\Sp^2$, such that $f_k\circ\phi_k$ is uniformly Lipschitz bounded and converges uniformly to $f$.
\end{thm}

{\small
Ernst Kuwert, Mathematisches Institut, Universit\"at Freiburg, Eckerstra{\ss}e 1, 
79104 Freiburg, Germany\\
{\sc Email: } ernst.kuwert@math.uni-freiburg.de\\
\\ 
Andrea Mondino, Scuola Normale Superiore, Piazza dei Cavalieri 7, 56126 Pisa, Italy\\
{\sc Email: } andrea.mondino@sns.it\\
\\
Johannes Schygulla, Mathematisches Institut, Universit\"at Freiburg, Eckerstra{\ss}e 1, 
79104 Freiburg, Germany\\
{\sc Email: } johannes.schygulla@math.uni-freiburg.de\\
}

\end{document}